\documentclass[10pt,twoside]{amsart}
\usepackage{lmodern}
\usepackage[T1]{fontenc}
\usepackage{amsmath,amsfonts,amsthm,amssymb,amscd}
\usepackage{mathrsfs}
\usepackage{latexsym}
\usepackage[pdftex]{graphicx}
\usepackage{enumitem}
\usepackage{bbm}
\usepackage{accents}

\usepackage{geometry}
\newtheorem{corollary}{Corollary}[section]
\newtheorem{lemma}{Lemma}[section]
\newtheorem{theorem}{Theorem}[section]

\usepackage{hyperref}  %added
\hypersetup{colorlinks=true, linktoc=all,
    linkcolor=blue,
		citecolor=blue,
		urlcolor=cyan,
		bookmarksopen=true
		}
\usepackage{hyperref}  %added
\newcommand{\lnc}{\mathscr{L}}
\def\T{\mathbb{T}}

\title{Waring's problem with almost proportional summands}
\author{Zarullo Rakhmonov, Firuz Rakhmonov}
\address{A.Dzhuraev Institute of Mathematics,  National Academy of Sciences of Tajikistan}
\email{zarullo.rakhmonov@gmail.com, rakhmonov.firuz@gmail.com}
\date{}

\begin{document}

\begin{abstract}
For $n \geq 3$, an asymptotic formula is derived for the number of representations of a sufficiently large natural number $N$ as a sum of $r = 2^n + 1$ summands, each of which is an $n$-th power of natural numbers $x_i$, $i = \overline{1, r}$, satisfying the conditions
$$
 |x_i^n-\mu_iN|\le H,\qquad  H\ge N^{1-\theta(n,r)+\varepsilon},\qquad \theta(n,r)=\frac2{(r+1)(n^2-n)},
$$
where $\mu_1, \ldots, \mu_r$ are positive fixed numbers, and $\mu_1 + \ldots + \mu_n = 1$. This result strengthens the theorem of E.M.~Wright.

Bibliography: 22 references.
\end{abstract}

\maketitle

\section{Introduction}
Waring problem with almost proportional summands was first investigated by M.~E.~Wright~\cite{Wright-1933, Wright-1934}. For the number of representations of a sufficiently large number $N$ in the form
 \begin{equation}
\label{Formula-x1n+...+xrn=N}
x_1^n+x_2^n+\ldots+x_r^n=N,
\end{equation}
where $x_1, x_2, \dots, x_r$ are natural numbers, and
$$
|x_i^n-\mu_iN|\le N^{1-\theta},\qquad i=1,2,\ldots,r,\qquad \mu_1+\mu_2+\ldots+\mu_r=1,
$$
$\mu_1,\ldots,\mu_r$~---~positive fixed numbers, and the number $\theta=\theta(n,r)$ is determined from the relation
\begin{equation}\label{formula teta(n,r) Daemena}
\theta(n,r)=\frac1n\min\left(\frac{(r-2^n)(2^{n-1}+1)}{(nr+n-2^n-3)2^{n-1}+r},\frac{r-(n-2)2^{n-1}-4}{r+2^{n-1}-4}, \frac{r-2^{n-1}}{nr-2^{n-1}+n-1}\right),
\end{equation}
he found an asymptotic formula for
\begin{equation}\label{formula r(n) Daemena}
r\ge r_n=(n-2)2^{n-1}+5.
\end{equation}
From this, particularly for $n=3,\ 4,\ 5,\ 6,\ 7,\ 8,\ 9,\ 10$ with $r=(n-2)2^{n-1}+5$, we have the following results in the table:
\begin{center}
\begin{tabular}
{  |c             ||   c  |  c | c  | c   |  c  | c   | c  | c   | c  }\hline
   $n$             &  3   &  4 & 5  &  6  &  7  &  8  &  9 &  10  \\ \hline
$r=(n-2)2^{n-1}+5$ &  9   & 21 & 53 & 133 &325  & 773 &1797&4101  \\  \hline\vspace{2pt}
$\theta(n,r)$      &$\frac1{51}$ & $\frac1{100}$ & $\frac1{325}$ &$\frac1{966}$&$\frac1{2695}$&$\frac1{6279}$&$\frac1{18441}$&$\frac1{46090}$ \\ \hline
  \end{tabular}\vspace{4pt}\\
  Table 1.
\end{center}
M.~E.~Wright's theorem on the asymptotic formula in Waring's problem with almost proportional summands for $\mu_1=\ldots=\mu_r=\frac1r$, that is, for
\begin{equation}\label{formula usl PRS v prob Waringa}
\left(\frac{N}r-N^{1-\theta}\right)^\frac1n\le x_j\le\left(\frac{N}r-N^{1-\theta}\right)^\frac1n
\end{equation}
turns into a theorem on the asymptotic formula in Waring's problem with almost equal summands.
Note that the inequality (\ref{formula usl PRS v prob Waringa}) can be represented as
$$
N_1-H\le x_j\le N_1+H,
$$
where\begin{align*}
&N_1=\frac{\left(\frac{N}r+N^{1-\theta}\right)^\frac1n+\left(\frac{N}r-N^{1-\theta}\right)^\frac1n}2= \left(\frac{N}r\right)^\frac1n\left(1+O\left(N^{-2\theta}\right)\right),\\
&H=\frac{\left(\frac{N}r+N^{1-\theta}\right)^\frac1n-\left(\frac{N}r-N^{1-\theta}\right)^\frac1n}2= \frac{r^{1-\frac1n}}nN^{\frac1n-\theta}\left(1+O\left(N^{-2\theta}\right)\right).
\end{align*}
R.~Vaughan~\cite{Vaughan-1981}, using the method of van der Corput in combination with Hua Lo-Keng's estimate for complete exponential sums of the form,
$$
S_b(a,q)=\sum_{k=1}^qe\left(\frac{ak^n+bk}{q}\right),\qquad S(a,q)=S_0(a,q),
$$
for the sums of G.~Weyl in major arcs, proved:
\begin{align*}
T(\alpha,x)=\sum_{m\le x}e(\alpha m^n)=\frac{S(a,q)}{q}\int_0^xe\left(\lambda t^n\right)dt+O\left(q^{\frac 12+\varepsilon}\left(1+x^n|\lambda|\right)^{\frac 12}\right).
\end{align*}
Under the condition that $\alpha$ is very well approximated by a rational number with denominator $q$, that is, under the condition
\begin{align*}
|\lambda|\leq \frac{1}{2nqx^{n-1}},
\end{align*}
he also proved:
\begin{align}\label{formula Vauhgan dlya T(alpha,x)}
&T(\alpha ,x)=\frac{x \ S(a,q)}{q}\int_0^1e\left(\lambda t^n\right)dt+O\left(q^{\frac 12+\varepsilon}\right).
\end{align}
R.~Vaughan's theorem for the case $n=3$ was generalized~\cite{Rakh-2008-1,Rakh-2014} in 2008~for short exponential sums of G.~Weyl of the form
\begin{align*}
&T(\alpha;x,y)=\sum_{x-y<m\le x}e(\alpha m^n),
\end{align*}
and, using this result, an asymptotic formula ~\cite{Rakh-2008-2,Rakh-2008-3,Rakh-2009} was established for the number of representations of a sufficiently large number $N$ in the form (\ref{Formula-x1n+...+xrn=N}) for $n=3$ and $r=9$, under the conditions
$$
\left|x_i-\left(\frac N9\right)^\frac13\right|\le H, \qquad i=1,2,\ldots,9, \qquad H\ge N^{\theta(3,9)+\varepsilon},\qquad \theta(3,9)=\frac3{10}.
$$
Then, R.~Vaughan's theorem was generalized for an arbitrary fixed $n$~\cite{Rakh-2010,Rakh-2015-1}. Using this result, an asymptotic formula was found for the number of solutions of the Diophantine equation (\ref{Formula-x1n+...+xrn=N}) for $n= 4$ and $n=5$ ~\cite{Rakh-2011,Rakh-2015-2,Rakh-2015-1}, under the conditions
$$
\left|x_i-\left(\frac{N}{2^n+1}\right)^{\frac{1}{n}}\right|\le H,\qquad i=1,\ldots,r,\qquad r=r(n)=2^n+1,\qquad H\ge N^{\frac{1}{n}-\theta(n,r)+\varepsilon},
$$
where
$$
\theta(4,17)=\frac{1}{108},\qquad \theta(5,33)=\frac{1}{340}.
$$
These results are a strengthening of E.~M.~Wright's theorem for $n=3,4,5$.

In 2010, D. Daemen \cite{Daemen-2010-1} using the mean value theorem of I. M. Vinogradov and the <<binomial descent>> procedure, for $r\ge r_n$, where
\begin{equation}\label{formula rez Daemena}
\begin{split}
&r_2=9,\quad  r_3=19,\quad  r_4=49,\quad  r_5=113,\quad r_6=243,\quad  r_7=417, \quad  r_8=675,\quad  r_9=1083,\\
&r_{10}=1773, \qquad r_n=2\left[\frac{5n^2}3\ln n+\frac{29n^2}{30}\ln\ln n+\frac{7n^2}3\ln\ln\ln n+Cn^2\right]+1,\quad \ n>10,
\end{split}
\end{equation}
$C$~---~absolute constant, proved an asymptotic formula for the number of solutions of the equation (\ref{Formula-x1n+...+xrn=N}), under the conditions
$$
X-Y\le x_j\le X+Y,\quad  1\le j\le r,\ \ X=\left[\left(\dfrac{N}r\right)^\frac1n\right],\quad  Y=\sqrt{X}Y_n,\quad Y_n=\left(\ln X\right)^{r_{n-1}},
$$
where $Y_2$ is a function of $X$ tending to infinity together with $X$ (see also \cite{Daemen-2010-2,Daemen-2010-3}). This asymptotic formula in the sense of the number of terms is weaker than the results of Z.~Kh.~Rakhmonov for $n=3,\ 4,\ 5$ and the results of E.~M.~Wright for $n= 6,\ 7$. It is stronger than the theorem of E.~M.~Wright for $n\ge8$.

In this work, it is proven that E. M. Wright's theorem on the asymptotic formula in the generalization of the Waring problem with almost proportional summands holds under the condition
\begin{equation}\label{formula paramenri teta i r}
\theta(n,r)=\frac2{(r+1)(n^2-n)}, \qquad r=2^n+1.
\end{equation}\vspace{2pt}

{\theorem\label{TeorAsForWaringPPSl} Let $N$ be a sufficiently large natural number, $n\ge3$ be a natural number, $r=2^n+1$, and $\mu_1,\ldots,\mu_r$ be positive fixed numbers satisfying the condition
$$
\mu_1+\ldots+\mu_r=1,
$$
$J_{n,r}(N,H)$ is the number of solutions to the Diophantine equation (\ref{Formula-x1n+...+xrn=N}) under the conditions
\begin{equation}\label{formula |xin-N/n|<=H}
 |x_i^n-\mu_iN|\le H,\qquad i=1,\ldots,r\qquad\theta(n,r)=\frac2{(r+1)(n^2-n)}.
\end{equation}
Then, for $H\ge N^{1-\theta(n,r)+\varepsilon}$, the following asymptotic formula holds:
$$
J_{n,r}(N,H)=\frac{2^r\gamma(n,r)}{n^r}\prod_{i=1}^r\mu_i^{-1+\frac1n}\mathfrak{S}(N)\frac{H^{r-1}}{N^{r-\frac rn}}+ O\left(\frac{H^{r-1}}{N^{r-\frac rn}\lnc}\right),
$$
where $\gamma(n,r)$ is an absolute constant determined by the relation
\begin{align*}
\gamma(n,r)=&\frac{r^{r-1}-\frac{r}{1!}(r-2)^{r-1}+\frac{r(r-1)}{2!}(r-4)^{r-1}-\frac{r(r-1)(r-2)}{3!}(r-6)^{r-1}+\ldots}{2^r(r-1)!},
\end{align*}
$\mathfrak{S}(N)$ is a singular series, the sum of which exceeds some positive constant, and the constant under the $O$ term depends on the values of $\mu_1,\ldots,\mu_r$.}

From this theorem, in particular, we obtain
\begin{center}
\begin{tabular}
{  |c     ||   c  |  c | c  | c   |  c  | c   | c  | c   | c  }\hline
   $n$           &  3 &  4 & 5  &  6  &  7  &  8  &  9 &  10  \\ \hline
$r=2^n+1$ &  9   & 17 & 33 & 64 & 129 & 257 &513&1025  \\  \hline\vspace{2pt}
$\theta(n,r)$    &$\frac1{30}$ & $\frac1{108}$ & $\frac1{340}$ &$\frac1{990}$&$\frac1{2730}$&$\frac1{7224}$&$\frac1{18504}$&$\frac1{46170}$\\ \hline
  \end{tabular}\vspace{4pt}\\
   Table 2.
\end{center}

From Theorem \ref{TeorAsForWaringPPSl}, an asymptotic formula for the generalization of the Waring problem with almost equal summands follows.
{\corollary \label{TeorAsForEstKubPr} Let $N$ be a sufficiently large natural number, $n\ge3$ be a natural number, $r=2^n+1$, and $J_{n,r}(N,H)$ be the number of solutions to the Diophantine equation (\ref{Formula-x1n+...+xrn=N}) under the conditions
$$
 \left|x_i^n-\frac{N}r\right|\le H,\qquad i=1,\ldots,r\qquad\theta(n,r)=\frac2{(r+1)(n^2-n)}.
$$
Then, for $H\ge N^{1-\theta(n,r)+\varepsilon}$, the following asymptotic formula holds:
$$
J_{n,r}(N,H)=\frac{2^rr^{r-\frac rn}\gamma(n,r)}{n^r}\mathfrak{S}(N)\frac{H^{r-1}}{N^{r-\frac rn}}+ O\left(\frac{H^{r-1}}{N^{r-\frac rn}\lnc}\right).
$$}

Note that theorem \ref{TeorAsForWaringPPSl} is a strengthening of the theorem of E.\,M.~Wright, and from the formula (\ref{formula rez Daemena}) and table 2 it also follows that the corollary \ref{TeorAsForEstKubPr} is stronger than the Dirk Damon theorem in the sense of the number of terms at least for $n=3,\ 4,\ 5,\ 6,\ 7, \ 8, \ 9, \ 10$.

The proof of the theorem \ref{TeorAsForEstKubPr} is carried out by the Hardy-Littlewood-Ramanujan circle method in the form of exponential sums of I.~M.~Vinogradov and for $n=3$ was previously proved in the paper~\cite{Rakh-2023}. The main statements that allowed obtaining new values (\ref{formula paramenri teta i r}) for the parameters $\theta(n,r)$ and $r$ include:
\begin{itemize}\vspace{4pt}
\item an asymptotic formula for short exponential sums of Weyl type $T(\alpha;x,y)$ in the small neighborhood of the center of major arcs (Corollary \ref{Sledst1 Teor ob poved kor trig summi Weyl} of Theorem \ref{Teor ob poved kor trig summi Weyl});\vspace{4pt}
\item a nontrivial estimate for sums $T(\alpha;x,y)$ in major arcs excluding the small neighborhood of their centers (Corollary \ref{Sledst2 Teor ob poved kor trig summi Weyl} of Theorem \ref{Teor ob poved kor trig summi Weyl});\vspace{4pt}
\item a nontrivial estimate for sums $T(\alpha;x,y)$ in minor arcs (Lemma \ref{Lemma-Weyl-3});\vspace{4pt}
\item  Theorem \ref{Teorema-obobsh-Hua} on the average value of short exponential sums of Weyl type $T(\alpha;x,y)$.\vspace{4pt}
\end{itemize}\vspace{2pt}
{\theorem \label{Teor ob poved kor trig summi Weyl} Let $\tau\geq 2n(n-1)x^{n-2}y$ and $\lambda\ge 0$, then for ${n|\lambda|x^{n-1}}\le \frac{1}{2q}$, the formula
\begin{equation*}
T(\alpha,x,y)=\frac{S(a,q)}{q}T(\lambda;x,y)+O\left(q^{\frac{1}{2}+\varepsilon}\right),
\end{equation*}
holds, and for $\{ n|\lambda|x^{n-1} \}>\frac{1}{2q}$, the following estimate holds
\begin{align*}
|T(\alpha,x,y)|&\ll q^{1-\frac 1n}\ln q +\min(yq^{-\frac 1n}, \lambda^{-\frac12} x^{1-\frac n2}q^{-\frac1n}).
\end{align*}}
{\corollary \label{Sledst1 Teor ob poved kor trig summi Weyl} Let $\tau\geq 2n(n-1)x^{n-2}y$, $|\lambda|\le \frac{1}{2nq x^{n-1}}$, then the next relation holds:
\begin{align*}
T(\alpha,x,y)=\frac{y}{q}S(a,q)\gamma(\lambda;x,y)+O(q^{\frac 12+\varepsilon
}),\\
\gamma(\lambda;x,y)=\int_{-0,5}^{0,5}e\left(\lambda\left(x-\frac y2+yt\right)^n\right)dt.
\end{align*}}
\vspace{2pt}
{\corollary \label{Sledst2 Teor ob poved kor trig summi Weyl} Let $\tau\geq2n(n-1)x^{n-2}y$, $\frac{1}{2nq x^{n-1}}<|\lambda|\le\frac1{q\tau}$, then the following estimate holds:
$$
T(\alpha,x,y)\ll q^{1-\frac1n}\ln q +\min\left(yq^{-\frac1n}, x^\frac12q^{\frac12-\frac1n}\right).
$$}
Corollary \ref{Sledst1 Teor ob poved kor trig summi Weyl} is a generalization of the formula (\ref{formula Vauhgan dlya T(alpha,x)}) for short Weyl exponential sums of the form $T(\alpha;x,y)$. The special case of Theorem \ref{Teor ob poved kor trig summi Weyl} for $n=3$ was previously proven in \cite{Rakh-2014} and is a refinement of Theorem 1 in \cite{Rakh-2015-2}. The proof of Theorem \ref{Teor ob poved kor trig summi Weyl} is carried out using the method of estimating special exponential sums by Van der Corput, applying the Poisson summation formula~\cite{Karatsuba}, and using the estimate by Hua Loo-Keng for complete exponential sums $S_b(a,q)$ (Lemma \ref{Lemma otsenka Hua o polnoy rats trigsummi}).

{\theorem\label{Teorema-obobsh-Hua} Let $x$ and $y$ be natural numbers, with $\sqrt{x} < y \leq x {\mathscr L}^{-1}$. Then, the estimate
 $$
 \int_{0}^{1}\left|T(\alpha;x,y)\right|^{2^k}d\alpha \ll y^{2^k-k+\varepsilon}, \qquad 1\le k\le n,
 $$
 holds.}

This theorem generalizes the result of Hua Loo-Keng~(\cite{Vaughan-1985}, lemma 2.5) on the average value of the exponential sum $T(\alpha,x)$, i.e., the estimate
$$
\int_0^1|T(\alpha,x)|^{2^k}\ll x^{2^k-k+\varepsilon},\qquad 1\le k\le n.
$$
\section{Well-known lemmas}
{\lemma \label{Lemma formula summirov Puassona}  {\rm \cite{Karatsuba}}. Let $f(u)$ be a real function such that $f''(u) > 0$ in the interval $[a, b]$, and let $\alpha$, $\beta$, and $\eta$ be arbitrary numbers satisfying $\alpha \le f'(a) \le f'(b) \le \beta$ and $0 < \eta < 1$. Then,
$$
\sum_{a<n\le b}e(f(n))=\sum_{\alpha-\eta<h\le\beta+\eta}\int_a^be(f(u)-hu)du+O(\eta^{-1}+\ln(\beta-\alpha+2)),
$$
where the constant in the $O$ notation is absolute.}
\vspace{2pt}
{\lemma \label{Lemma otsenka Hua o polnoy rats trigsummi} {\rm \cite{Hua}}. Let $(a, q) = 1$, where $q$ is a natural number, $b$ is any integer, and $\delta$ is any positive number not exceeding $0.00001$. Then, we have
$$
 S_b(a,q)=\sum_{k=1}^qe\left(\frac{ak^n+bk}{q}\right)\ll q^{\frac12+\delta}(b,q).
 $$}\vspace{4pt}

 {\lemma\label{Lemma otsenka trig integ po pervoy proizv}{\rm \cite{Arkhipov}}. Let the real function $f(u)$ and the monotonic function $g(u)$ satisfy the conditions: $f'(u)$ is monotonic, $|f'(u)| \geq m_1 > 0$, and $|g(u)| \leq M$. Then, the estimate holds:
$$
\int_a^bg(u)e(f(u))du\ll \frac{M}{m_1}.
$$}\vspace{2pt}

{\lemma\label{Lemma otsenka polnoy trigsumm}{\rm \cite{Arkhipov}}. Let $(a, q) = 1$, where $q$ is a natural number. Then, we have
$$
S(a,q)=\sum_{k=1}^{q}e\left(\frac{ak^n}{q}\right)\ll q^{1-\frac1n},
$$
where the constant under the Vinogradov symbol depends on $n$.}\vspace{2pt}
{\lemma \label{Lemma otsenka trig int po n-oy proizvodn}{\rm \cite{Arkhipov}}. Let $f(u)$ be a real function with a $n$-th order derivative ($n>1$) in the interval $a\le u\le b$. Suppose there exists a constant $A>0$ such that the inequality $A\le |f^{(n)}(u)|$ holds.
Then, the following estimate holds:
$$
\int\limits_a^be(f(u))du\le \min (b-a,6nA^{-\frac 1n}).
$$}
{\lemma\label{Lemma tochnoe znach trig integ}{\rm (\cite{Whittaker} page 174.)} For $m > 0$ and a natural number $r$, the following formula holds:
\begin{align*}
\int\limits_0^\infty\frac{\sin^rmt}{t^r}dt=&\frac{\pi m^{m-1}}{2^r(r-1)!}\left(r^{r-1}-\frac{r}{1!}(r-2)^{r-1}+\right.\\
 &\left.+\frac{r(r-1)}{2!}(r-4)^{r-1}-\frac{r(r-1)(r-2)}{3!}(r-6)^{r-1}+\ldots\right).
\end{align*}}
{\lemma\label{Lemma-Weyl-3} {\rm \cite{Rakh-2018}}. Let $x\ge x_0>0$, $y_0<y\le 0.01x$, $\tau(\cdot)$~---~the divisor function, $\alpha$~---~a real number,
$$
\left|\alpha-\frac{a}{q}\right|\le \frac{1}{q^2}, \qquad (a,q)=1,
$$
then the following estimate holds:
\begin{align*}
|T(\alpha;x,y)|\le 2y\left(4n!\left(\frac{1}q+\frac{1}y+\frac{q\ln q}{y^n}\right)\max_{h<y^{n-1}}\tau(h)\right)^{\frac{1}{2^{n-1}}},
\end{align*}}
This estimate is nontrivial when $q\gg 2^{2n-1}4n!\tau(y^{n-1})$, that is, when $q\gg y^\varepsilon$.
\vspace{4pt}

\section{
Proof of Theorem \ref{Teor ob poved kor trig summi Weyl} on the\\
behavior of short Weyl exponential sums in major arcs}
Using the orthogonality property of the complete linear rational exponential sum, we obtain
\begin{align}
T(\alpha;x,y)&=\sum_{x-y<m\leq x}e\left(\frac{ak^n}{q}+\lambda m^n\right)\sum_{\substack{k=0\\ k\equiv m\bmod q}}^{q-1}1
=\frac{1}{q}\sum_{b=0}^{q-1}T_b(\lambda;x,y)S_b(a,q), \label{formula T(alpha;x,y)=1/qsummiTb(alpha;x,y)Sb(a,q)}
 \end{align}
where
\begin{align*}
T_b(\lambda;x,y)=\sum_{x-y<m\leq x}e\left(\lambda m^n-\frac{bm}{q}\right), \qquad T(\lambda;x,y)=T_0(\lambda;x,y).
\end{align*}
Furthermore, without loss of generality, we will assume that $\lambda \ge 0$. The case $\lambda\le 0$ can be reduced to the case $\lambda\ge0$ by applying the formula (\ref{formula T(alpha;x,y)=1/qsummiTb(alpha;x,y)Sb(a,q)}) along with the transformation:
\begin{align*}
\overline{T(\alpha;x,y)}=\frac{1}{q}\sum_{b=0}^{q-1}T_{q-b}(-\lambda;x,y)S_{q-b}(q-a,q)= \frac{1}{q}\sum_{b=0}^{q-1}T_b(-\lambda;x,y)S_b(q-a,q).
\end{align*}
Keeping in mind that $n\lambda x^{n-1}- \{n\lambda x^{n-1} \}$ is an integer, let's express $T_b(\lambda;x,y)$ as follows:
\begin{align*}
&T_b(\lambda;x,y)=\sum_{x-y<m\le x}e\left(f_b(m)\right), \\
&f_b(m)=\lambda m^n-(n\lambda x^{n-1}-\{n\lambda x^{n-1}\})m-\frac{bm}q.
\end{align*}
We find the first and second order derivatives of the function $f_b(m)$:
$$
f_b'(m)=n\lambda(m^{n-1}-x^{n-1})+\{n\lambda x^{n-1}\}-\frac bq,\qquad  f_b''(m)=n(n-1)\lambda m^{n-2}\ge0.
$$
Therefore, the function $f_b'(m)$, $m\in (x-y,x]$ is non-decreasing. Thus, for any $b$, $b=0,1,\ldots,q-1$, the following inequality holds:
\begin{equation}\label{formula f'(x-y,b)<f_b'(m)<=f'(x,b)}
f_b'(x-y)<f_b'(m)\le f_b'(x).
\end{equation}
Estimating $f_b'(x)$ from above, we have:
\begin{equation}\label{formula f'(x,b)<1-}
f_b'(x)=\{n\lambda x^{n-1}\}-\frac bq<1-\frac{b}{q}.
\end{equation}
To estimate $f_b'(x-y)$ from below, using the inequality:
\begin{equation}\label{formula W=summ..}
W=\sum_{k=2}^{n-1}(-1)^{k}C_{n-1}^{k}x^{n-1-k}y^{k}\ge 0,\quad n\ge 3, \quad 3x\ge (n-3)y,
\end{equation}
and the conditions $|\lambda|\leq \frac{1}{q\tau}$ and $\tau \ge 2n(n-1)x^{n-2}y$, we have
\begin{align*}
f_b'(x-y)&=-n\lambda\left(x^{n-1}-\left(x-y\right)^{n-1}\right)+\{n\lambda x^{n-1}\}-\frac bq= \\
&=-n(n-1)\lambda x^{n-2}y+n\lambda W+\{n\lambda x^{n-1}\}-\frac bq\ge\\
&\ge-n(n-1)\lambda x^{n-2}y-\frac bq\ge-\frac{n(n-1)x^{n-2}y}{q\tau}-\frac{b}{q}\ge-1+\frac{1}{2q}.
\end{align*}
From here, using (\ref{formula f'(x,b)<1-})  and (\ref{formula f'(x-y,b)<f_b'(m)<=f'(x,b)}), we obtain
$$
-1+\frac1{2q}<f_b'(m)<1-\frac1q.
$$
Taking this inequality into account and applying the Poisson summation formula (lemma \ref{Lemma formula summirov Puassona}) to the sum $T_{b}(\lambda;x,y)$ with $\alpha =-1$, $\beta =1$, $\eta=0,5$, we obtain
\begin{align}\label{formula Tb(lambda;x,y)=Ib(-1)+Ib(0)+Ib(1)+O(1)}
&T_b(\lambda;x,y)=\sum_{h=-1}^1I_b(h)+O(1),\\
I_b(h)&=\int_{x-y}^xe(f_b(u,h))du,\quad  f_b(u,h)=f_b(u)-hu. \nonumber
\end{align}
The function
$$
f'_b(u,h)=n\lambda (u^{n-1}-x^{n-1})+\{n\lambda x^{n-1}\}-\frac{b}{q}-h
$$
on the interval $u\in [x-y,x]$ is a non-decreasing function. Therefore, the next inequality holds:
$$
f'_b(x-y,h)\le f'_b(u,h)\le f'_b(x,h).
$$
Using the formula (\ref{formula W=summ..}), as well as the conditions $|\lambda|\leq \frac1{q\tau}$ and $\tau\ge2n(n-1)x^{n-2}y$, this inequality can be expressed as follows:
\begin{align}
\label{otsenka f'h(u,b)}
& \{n\lambda x^{n-1}\}-\frac{b}{q}-h-\eta<f'_b(u,h)\le \{n\lambda x^{n-1}\}-\frac{b}{q}-h,\\
&\eta=n(n-1)\lambda x^{n-2}y-n\lambda W\le n(n-1)\lambda x^{n-2}y\le\frac{n(n-1)x^{n-2}y}{q\tau}\le\frac1{2q}.\nonumber
\end{align}
Furthermore, substituting (\ref{formula Tb(lambda;x,y)=Ib(-1)+Ib(0)+Ib(1)+O(1)}) for $b\neq0$ into (\ref{formula T(alpha;x,y)=1/qsummiTb(alpha;x,y)Sb(a,q)}), we obtain:
\begin{align}
\label{formula T(alpha;x,y)=S0(a,q)T0(lambda;x,y)/q+T-1+T0+T1+R}
T(\alpha;x,y)=\frac{S_0(a,q)}{q}T_0(\lambda;x,y)+\sum_{h=-1}^1\mathcal{T}(h)+\mathcal{R},\\
\mathcal{T}(h)=\frac{1}{q}\sum_{b=1}^{q-1}I_b(h)S_b(a,q),\qquad \mathcal{R}\ll\frac{1}{q}\sum_{b=1}^{q-1}|S_b(a,q)|. \nonumber
\end{align}
Using lemma \ref{Lemma otsenka Hua o polnoy rats trigsummi}  with $\delta=0,5\varepsilon$, let's estimate the remainder term $\mathcal{R}$:
\begin{align}
\mathcal{R}\le&\frac{1}{q}\sum_{b=1}^{q-1}|S_b(a,q)|\ll q^{-\frac12+\delta}\sum_{b=1}^{q-1}(b,q)=q^{-\frac12+\delta}\sum_{d\backslash q}d\sum_{\substack{1\le b\le q-1\\(b,q)=d}}1=\nonumber\\
&=q^{-\frac 12+\delta}\sum_{d\backslash q}d\varphi\left(\frac qd\right)\le q^{\frac12+\delta}\tau(q)\ll q^{\frac12+\varepsilon}.\label{formula otsenka R<<q^1/2+varepsilon}
\end{align}
Estimating from above the sums $\mathcal{T}(1)$ and $\mathcal{T}(-1)$, and setting $h=1$ in (\ref{otsenka f'h(u,b)}), we obtain:
\begin{align*}
&f'_b(u,1)\le\{n\lambda x^{n-1}\}-\frac{b}{q}-1\le -\frac{b}{q}< 0,
\end{align*}
and estimating the integral $|I_b(1)|$ in terms of the first derivative (lemma \ref{Lemma otsenka trig integ po pervoy proizv}), we have:
$$
|I_b(1)|=\left|\int_{x-y}^xe(f_b(u,1))du\right|\ll \frac{q}{b}.
$$
Using this estimate and then lemma \ref{Lemma otsenka Hua o polnoy rats trigsummi} with $\delta=0.5\varepsilon$, we have:
\begin{align}
\mathcal{T}(1)&=\frac{1}{q}\sum_{b=1}^{q-1}I_b(1)S_b(a,q)\ll\sum_{b=1}^{q-1}\frac{|S_b(a,q)|}{b}\ll
 q^{\frac12+\delta}\sum_{b=1}^{q-1}\frac{(b,q)}{b}\ll q^{\frac{1}{2}+\varepsilon}.\label{formula otsenka T1<<q^1/2+varepsilon}
\end{align}
Substituting $h=-1$ into (\ref{otsenka f'h(u,b)}), we have:
\begin{align*}
f'_b(u,-1)&>\{n\lambda x^{n-1}\}+\frac{q-b}{q}-\eta\ge \frac{q-b}{q}.
\end{align*}
Let's estimate the integral $I(-1,b)$ using the magnitude of the first derivative (lemma \ref{Lemma otsenka trig integ po pervoy proizv}). We obtain:
$$
|I_b(-1)|=\left|\int_{x-y}^xe(f_b(u,-1))du\right|\ll\frac{q}{q-b}.
$$
From here, proceeding similarly to the case of estimating $\mathcal{T}(1)$, we have:
\begin{align}
\mathcal{T}(-1)=&\sum_{b=1}^{q-1}\frac{I_b(-1)S_b(a,q)}{q}\ll\sum_{b=1}^{q-1}\frac{|S_b(a,q)|}{q-b}\ll q^{\frac12+\delta} \sum_{b=1}^{q-1}\frac{(b,q)}{b}\ll q^{\frac12+\varepsilon}.\label{formula otsenka T-1<<q^1/2+varepsilon}
\end{align}
Substituting the estimates for $\mathcal{T}(1)$, $\mathcal{T}(-1)$, and $\mathcal{R}$ from (\ref{formula otsenka T1<<q^1/2+varepsilon}), (\ref{formula otsenka T-1<<q^1/2+varepsilon}), and (\ref{formula otsenka R<<q^1/2+varepsilon}) respectively into the right-hand side of (\ref{formula T(alpha;x,y)=S0(a,q)T0(lambda;x,y)/q+T-1+T0+T1+R}), we obtain:
\begin{align}
\label{formula T(alpha;x,y)=S0(a,q)T0(lambda;x,y)/q+T0+..}
T(\alpha;x,y)=\frac{S_0(a,q)}{q}T_0(\lambda;x,y)+\mathcal{T}(0)+O\left(q^{\frac{1}{2}+\varepsilon}\right).
\end{align}
Now, using this relation, we will separately prove the first and second statements of the theorem, corresponding to the cases $\{n\lambda x^{n-1}\}\leq\frac{1}{2q}$ and $\{n\lambda x^{n-1} \}>\frac{1}{2q}$.

\textbf{1. Case $ \{n\lambda x^{n-1} \}\le \frac{1}{2q}$.} Let's estimate $\mathcal{T}(0)$. To do this, setting $h=0$ in (\ref{otsenka f'h(u,b)}), we have
\begin{align*}
&f'_b(u,0)\le\{n\lambda x^{n-1}\}-\frac{b}{q}\le\frac{1-2b}{2q}\le-\frac{b}{2q}<0.
\end{align*}
Using this inequality and estimating the integral $I_b(0)$ by the magnitude of its first derivative (lemma \ref{Lemma otsenka trig integ po pervoy proizv}), we find
$$
|I_b(0)|=\left|\int_{x-y}^xe(f_b(u,0))du\right|\ll\frac{q}{b}.
$$
Proceeding similarly to the case of estimating $\mathcal{T}(1)$, we obtain
$$
\mathcal{T}(0)=\sum_{b=1}^{q-1}\frac{I_b(0)S_b(a,q)}{q}\ll\sum_{b=1}^{q-1}\frac{|S_b(a,q)|}{b}\ll q^{\frac12+\delta} \sum_{b=1}^{q-1}\frac{(b,q)}{b}\ll q^{\frac12+\varepsilon}.
$$
Substituting this estimate into the right-hand side of (\ref{formula T(alpha;x,y)=S0(a,q)T0(lambda;x,y)/q+T0+..}), and also noting that
$$
S_0(a,q)=S(a,q)\qquad T(\lambda;x,y)=T_0(\lambda;x,y),
$$
we obtain the first statement of the theorem.

\textbf{2. Case $\{n\lambda x^{n-1}\}>\frac{1}{2q}$.} Applying formula (\ref{formula Tb(lambda;x,y)=Ib(-1)+Ib(0)+Ib(1)+O(1)}) to the sum $T_0(\lambda;x,y)$ in the relation (\ref{formula T(alpha;x,y)=S0(a,q)T0(lambda;x,y)/q+T0+..}), and then proceeding to estimates and using the estimate $S_b(a,q)\ll q^{1-\frac1n}$ (lemma \ref{Lemma otsenka polnoy trigsumm}), we find:
\begin{align}
T(\alpha;x,y)&=\frac{S_0(a,q)}{q}\left(\sum_{h=-1}^1I_0(h)+O(1)\right)+\mathcal{T}(0)+O\left(q^{\frac12+\varepsilon}\right)=\nonumber\\
&=\frac{S_0(a,q)}{q}\left(I_0(-1)+I_0(1)\right)+\sum_{b=0}^{q-1}\frac{I_b(0)S_b(a,q)}q+O\left(q^{\frac12+\varepsilon}\right)\ll\nonumber\\
&\ll q^{-\frac1n}\left(I_0(-1)+I_0(1)+\mathcal{R}(0)\right)+q^{\frac12+\varepsilon},\label{formula T(alpha;x,y)<<}\\
&\qquad\mathcal{R}(0)=\sum_{b=0}^{q-1}I_b(0).\label{formula R0=..}
\end{align}
In the relation (\ref{otsenka f'h(u,b)}), setting $b=0$ and $h=-1$, we have:
\begin{align*}
f'_0(u,-1)&>\{n\lambda x^{n-1}\}+1-\eta\ge1-\frac1{2q}\ge\frac12.
\end{align*}
Utilizing this inequality and estimating the integral $I_0(-1)$ in terms of the first derivative (lemma \ref{Lemma otsenka trig integ po pervoy proizv}), we obtain:
\begin{equation}\label{formula otsenka I_0(-1)<<1}
|I_0(-1)|=\left|\int_{x-y}^xe(f_0(u,-1))du\right|\ll 1.
\end{equation}
We will estimate the integral $I_0(1)$ in terms of the second-order derivative (lemma \ref{Lemma otsenka trig int po n-oy proizvodn}). To do this, using the inequality $f_0''(u,1)\ge n(n-1)\lambda (x-y)^{n-2}\gg \lambda x^{n-2}$, we find:
\begin{equation}\label{formula otsenka I_0(1)<<mim(y,..)}
|I_0(1)|=\left|\int_{x-y}^xe(f_0(u,1))du\right|\ll\min\left(y,\lambda^{-\frac12} x^{1-\frac n2}\right).
\end{equation}
Now let's estimate the sum $\mathcal{R}(0)$. We define a natural number $r$ by the relation:
\begin{equation*}
\frac{r}{2q}\le \{n\lambda x^{n-1}\}<\frac{r+1}{2q},  \qquad 1\le r\le 2q-1.
\end{equation*}
From this, using the inequality (\ref{otsenka f'h(u,b)}) with $h=0$ and the condition $\eta \le \frac{1}{2q}$, we obtain:
\begin{align}
&f'_b(u,0)>\{n\lambda x^{n-1}\}-\frac{b}{q}-\eta\ge \frac{r-2b-1}{2q},\label{otsenka f'b(u,0) snizu}\\
&f'_b(u,0)\le\{n\lambda x^{n-1}\}-\frac{b}{q}<\frac{r-2b+1}{2q}. \label{otsenka f'b(u,0) sverkhu}
\end{align}
Let $r=2\rho$ be even ($1\le\rho\le q-1$). We divide the range of summation $0\le b \le q-1$ in the sum $\mathcal{R}(0)$ into the following three sets:
$$
0\le b\le\rho-1, \qquad  b=\rho, \qquad \rho+1\le b\le q-1,
$$
where in the first set the right-hand side of inequality (\ref{otsenka f'b(u,0) snizu}) is greater than zero, and in the third set the right-hand side of inequality (\ref{otsenka f'b(u,0) sverkhu}) is less than zero. That is,
\begin{align*}
&f'_b(u,0)>\frac{2\rho-2b-1}{2q}\ge \frac{\rho-b}{2q}, &0\le b\le\rho-1,\\
&f'_b(u,0)<\frac{2\rho-2b+1}{2q}\le\frac{\rho-b}{2q}, &\rho+1\le b\le q-1.
\end{align*}
Using these inequalities and estimating the integral $I_b(0)$ in terms of the first derivative, we find
$$
I_b(0)=\int_{x-y}^xe(f_b(u,0))du\ll\frac{q}{|\rho-b|}, \quad b\neq \rho.
$$
In the case $b=\rho$, using the relation
$$
f_{\rho}''(u,)=n(n-1)\lambda(x-y)^{n-2}\gg\lambda x^{n-2},
$$
and estimating the integral $I_{\rho}(0)$ in terms of the second-order derivative (Lemma \ref{Lemma otsenka trig int po n-oy proizvodn}), we find
\vspace{-4pt}
$$
|I_{\rho}(0)|\ll\min\left(y,\lambda^{-\frac12} x^{1-\frac n2}\right).
$$
\vspace{-2pt}
Substituting the obtained estimates for $I_b(0)$ into (\ref{formula R0=..}), we get
\vspace{-4pt}
\begin{align*}
\mathcal{R}(0)&=\sum_{\substack{b=0\\ b\neq\rho}}^{q-1}I_b(0)+I_{\rho}(0)\ll\sum_{\substack{b=0\\ b\neq\rho}}^{q-1}\frac{q}{|\rho-b|}+\min\left(y,\lambda^{-\frac12}x^{1-\frac n2}\right)\ll q\ln q+\min\left(y,\lambda^{-\frac12}x^{1-\frac n2}\right).
\end{align*}\vspace{-4pt}

Now let $r=2\rho+1$ be odd ($0\le\rho\le q-1$). We split the summation interval $0\le b \le q-1$ in the sum $\mathcal{R}(0)$ into the following three sets:
\vspace{-3pt}
$$
0\le b\le\rho-1, \qquad  b=\rho,\, \rho+1, \qquad\rho+2\le b \le q-1,
$$
\vspace{-3pt}
in the first of which the right-hand side of the inequality (\ref{otsenka f'b(u,0) snizu}) is greater than zero, and in the third one the right-hand side of the inequality (\ref{otsenka f'b(u,0) sverkhu}) is less than zero. Thus, we have
\vspace{-3pt}
\begin{align*}
&f'_b(u,0)>\frac{2\rho-2b-1}{2q}\ge\frac{\rho-b}{2q}, &0\le b \le \rho-1,\\
&f'_b(u,0)<\frac{2\rho+1-2b+1}{2q}\le \frac{\rho-b}{2q}, &\rho+2\le b \le q-1.
\end{align*}
\vspace{-3pt}
Therefore,
\vspace{-3pt}
\begin{align*}
I_b(0)=\int_{x-y}^xe(f_b(u,0))du\ll\frac{q}{|\rho-b|}, \qquad b\neq\rho-1, \quad  b\neq\rho.
\end{align*}
\vspace{-3pt}
In the case where $b=\rho-1$ or $b=\rho$, following a similar approach to the previous estimate for $I_\rho(0)$, we find
$$
|I_b(0)|\ll\min\left(y,\lambda^{-\frac12} x^{1-\frac n2}\right).
$$
Substituting the obtained estimates for $I_b(0)$ into (\ref{formula R0=..}), we get
\vspace{-3pt}
\begin{align*}
\mathcal{R}(0)=&\sum_{\substack{b=0\\ b\neq\rho-1,\,\rho}}^{q-1}I_b(0)+I_{\rho-1}(0)+I_{\rho}(0)\ll\sum_{\substack{b=0\\ b\neq\rho-1,\,\rho}}^{q-1}\frac{q}{|\rho-b|}+\min\left(y,\lambda^{-\frac12}x^{1-\frac n2}\right)\ll\\
&\ll q\ln q+\min\left(y,\lambda^{-\frac12}x^{1-\frac n2}\right).
\end{align*}
\vspace{-3pt}
Substituting the estimates for $I_0(-1)$ and $I_0(1)$, respectively, from formulas (\ref{formula otsenka I_0(-1)<<1}) and (\ref{formula otsenka I_0(1)<<mim(y,..)}), and the estimate for $\mathcal{R}(0)$ into (\ref{formula T(alpha;x,y)<<}), we obtain the second statement of the theorem.
\vspace{-2pt}
\section{Proof of Theorem \ref{Teorema-obobsh-Hua} on the average value of the short exponential sum $T(\alpha;x,y)$}
{\lemma\label{Lemma razn oper} Let $\Delta_k$ denote the $k$-th application of the difference operator, so that for any real-valued function $f(u)$,
\begin{align}
&\Delta_1(f(u);h)=f(u+h)-f(u),\nonumber\\
&\Delta_{k+1}(f(u);h_1,\ldots,h_{k+1})=\Delta_1(\Delta_k(f(u);h_1,\ldots,t_k);h_{k+1}). \label{Opred Delta{k+1}(f(u);h1,...,h{k+1})}
\end{align}
Then for $k=1,\ldots,n-1$, the next relation holds
$$
\Delta_k(u^n;h_1,\ldots,h_k)=h_1\ldots h_kg_k(u;h_1,\ldots,h_k),
$$
where $g_k=g_k(u;h_1,\ldots,h_k)$ is a form of degree $n-k$ with integer coefficients, having a degree $n-k$ with respect to $u$ and the leading coefficient $n(n-1)\ldots(n-k+1)$, i.e.,
$$
g_k(u;h_1,\ldots,h_k)=\frac{n!}{(n-k)!}u^{n-k}+\ldots.
$$}
\begin{proof}
Applying the formula (\ref{Opred Delta{k+1}(f(u);h1,...,h{k+1})}), we find
\begin{align*}
&\Delta_1(u^n;h_1)=(u+h_1)^n-u^n=\sum_{i_1=1}^nC_n^{i_1}h_1^{i_1}u^{n-i_1},\\
&\Delta_2(u^n;h_1,h_2)=\Delta_1\left(\Delta_1(u^n;h_1);h_2\right)=\Delta_1\left(\sum_{i_1=1}^nC_n^{i_1}h_1^{i_1}u^{n-i_1};h_2\right)=\\
&=\sum_{{i_1}=1}^nC_n^{i_1}h_1^{i_1}(u+h_2)^{n-i_1}-\sum_{i_1=1}^nC_n^{i_1}h_1^{i_1}u^{n-i_1}= \sum_{{i_1}=1}^nC_n^{i_1}h_1^{i_1} \sum_{i_2=1}^{n-i_1}C_{n-i_1}^{i_2}h_2^{i_2}u^{n-i_1-i_2} .
\end{align*}
Sequentially applying the formula (\ref{Opred Delta{k+1}(f(u);h1,...,h{k+1})}), it can be easily shown, for $k=1,2,\ldots,n-1$, that the following equality holds:
\begin{align*}
&\Delta_k(u^n;h_1,\ldots,h_k)=\sum_{{i_1}=1}^nC_n^{i_1}h_1^{i_1}%\sum_{i_2=1}^{n-i_1}C_{n-i_1}^{i_2}h_2^{i_2}
\ldots \sum_{i_k=1}^{n-i_1-\ldots-i_{k-1}} C_{n-i_1-\ldots-i_{k-1}}^{i_k}h_k^{i_k}u^{n-i_1-\ldots-i_k}.
\end{align*}
From this formula, it follows that the next formula holds:
$$
\Delta_k(u^n;h_1,\ldots,h_k)=h_1\ldots h_kg_k(u;h_1,\ldots,h_k),
$$
where $g_k=g_k(u;h_1,\ldots,h_k)$ is a form of order $n-k$ with integer coefficients, having a degree $n-k$ with respect to $u$ and a leading coefficient of $n(n-1)\ldots(n-k+1)$, i.e.,
$$
g_k(u;h_1,\ldots,h_k)=\frac{n!}{(n-k)!}u^{n-k}+\ldots.
$$
\end{proof}
{\lemma\label{Lemma virazh st T(f(u);x,y) konech razn} Let $f=f(m)$ be a polynomial of degree $n$, $x$ and $y$ are positive integers with $y<x$,
$$
\T(f;x,y)=\sum_{x-y<m\le x}e(f(m)),
$$
then for $k=1,\ldots,n-1$, the following inequality holds:
\begin{align*}
&|\T;x,y)|^{2^k}\leq(2y)^{2^k-k-1}\sum_{|h_1|<y}\ldots\sum_{|h_k|<y}\T_k, \\
&\T_k=\left|\sum_{m\in I_k(x,y;h_1,\ldots, h_k)}e(\Delta_k(f(m); h_1,\ldots,h_k))\right|,
\end{align*}
where the intervals $I_k(x,y;h_1,\ldots,h_k)$ are defined by the relations:
\begin{align*}
&I_1(x,y;h_1)=(x-y,x]\cap (x-y-h_1,x-h_1],\\
&I_k(x,y;h_1,\ldots, h_{k})=I_{k-1}(x,y;h_1,\ldots, h_{k-1})\cap I_{k-1}(x-h_k,y;h_1,\ldots, h_{k-1}).
\end{align*}
That  is, the interval $I_{k-1}(x-h_k,y;h_1,\ldots, h_{k-1})$ is obtained from the interval \linebreak $I_{k-1}(x,y;h_1,\ldots, h_{k-1})$ by shifting all the intervals, of which it is an intersection, by an amount of $-h_k$.
}
\begin{proof} We proceed with the proof using mathematical induction on $k$. For $k=1$, we have
\begin{align*}
|\T(f;x,y)|^2&=\sum_{x-y<m\le x}\sum_{x-y-m<h\le x-m}\hspace{-20pt}e(f(m+h)-f(m))=\sum_{|h|<y}\sum_{m\in I_1(x,y;h_1)}e(\Delta_1(f(m);h))\le\sum_{|h|<y}\T_1.
\end{align*}
Assume that the lemma holds for $k$, $1\le k\le n-2$, i.e.,
$$
|\T(f(m);x,y)|^{2^{k}}\leq(2y)^{2^k-k-1}\sum_{|h_1|<y}\sum_{|h_2|<y}\ldots\sum_{|h_k|<y}\T_k.
$$
Squaring both sides of this inequality and then applying the Cauchy inequality successively to the sums over $h_1, \ldots, h_k$, we find
\begin{align}
|\T&(f(m);x,y)|^{2^{k+1}}\le(2y)^{2^{k+1}-(k+1)-1}\sum_{|h_1|<y}\ldots\sum_{|h_k|<y}\T_k^2. \label{aaa}
\end{align}
From the equivalence of the relations $m_1\in I_k(x,y;h_1,\ldots, h_{k})$ and $m_1-m\in I_k(x-m,y;h_1,\ldots, h_{k})$, we have
\begin{align*}
&\T_k^2=\sum_{\substack{m\in I_k(x,y;h_1,\ldots, h_{k})\\ m_1-m\in I_k(x-m,y;h_1,\ldots, h_{k})}}e(\Delta_k(f(m_1); h_1,\ldots,h_k)-\Delta_k(f(m); h_1,\ldots,h_k)).
\end{align*}
Denoting the difference $m_1-m$ as $h_{k+1}$ and then making the sum over $h_{k+1}$ outer, utilizing the equivalence of relations $h_{k+1}\in I_k(x-m,y;h_1,\ldots, h_{k})$ and $m\in I_k(x-h_{k+1},y;h_1,\ldots, h_{k})$, and using the relation (\ref{Opred Delta{k+1}(f(u);h1,...,h{k+1})}), we have
\begin{align*}
\T_k^2&=\sum_{|h_{k+1}|<y}\sum_{\substack{m\in I_k(x,y;h_1,\ldots, h_{k})\\m\in I_k(x-h_{k+1},y;h_1,\ldots, h_{k})}}e(\Delta_{k+1}(f(m);h_1,\ldots,h_{k+1}))=\\
&=\sum_{|h_{k+1}|<y}\sum_{m\in I_{k+1}(x,y;h_1,\ldots, h_{k+1})}e(\Delta_{k+1}(f(m);h_1,\ldots,h_{k+1}))\le\sum_{|h_{k+1}|<y}\T_{k+1}.
\end{align*}
Substituting the right-hand side of the last inequality into (\ref{aaa}), we obtain the statement of the lemma.
\end{proof}

{\sc Proof of Theorem \ref{Teorema-obobsh-Hua}}. We will use mathematical induction on $k$. For $k=1$, using the fact that for $x-y<m_1,m_2\le x$, the Diophantine equations $m_1^n=m_2^n$ and $m_1=m_2$ are equivalent, we have
\begin{align*}
\int_0^1&|T(\alpha;x,y)|^2 d\alpha=\sum_{x-y<m_1,m_2\le x}\int_0^1e(\alpha(m_1^n-m_2^n))d\alpha=%\sum_{\substack{x-y<m_1,m_2\le x\\ m_1^n=m_2^n}}1 =
\sum_{x-y<m_1\le x}1\ll y.
\end{align*}
Now, assume that the theorem holds for $2\le k\le n-1$, i.e.,
\begin{equation}\label{formula int|Tn(alpha;x,y)|2kdalpha<=y2k-k+varepsilon}
\int_0^1 |T(\alpha;x,y)|^{2^k}d\alpha \le y^{2k-k+\varepsilon}.
\end{equation}
In Lemma \ref{Lemma virazh st T(f(u);x,y) konech razn}, setting $f(m)=\alpha m^n$, we have
\begin{align*}
&|T(\alpha;x,y)|^{2^k}\leq(2y)^{2^k-k-1}\sum_{|h_1|<y}\ldots\sum_{|h_k|<y}\left|\sum_{m\in I_k}e\left(\alpha\Delta_k(m^n; h_1,\ldots,h_k)\right)\right|.
\end{align*}
Using Lemma \ref{Lemma razn oper}, we find
$$
\Delta_k(m^n;h_1,\ldots,h_k)=h_1\ldots h_kg_k(m;h_1,\ldots,h_k),
$$
where $g_k=g_k(m;h_1,\ldots,h_k)$ is a form of degree $n-k$ with integer coefficients, having a degree $n-k$ with a leading coefficient of $n(n-1)\ldots(n-k+1)$, i.e.,
$$
g_k(m;h_1,\ldots,h_k)=\frac{n!}{(n-k)!}m^{n-k}+\ldots.
$$
From this and the conditions $x-y<m\le x$, $|h_i|<y$, $i=1,\ldots,k$, $\sqrt{x}<y\le x\lnc^{-1}$, it follows that there exists $x_0$ such that for $x>x_0$, the next inequality holds:
\begin{equation}\label{formula gk(u;h1,...hk)>0}
g_k(m;h_1,\ldots,h_k)>0.
\end{equation}
Denoting by $r(h)$ the number of solutions to the Diophantine equation
$$
h_1\ldots h_kg_k(m;h_1,\ldots,h_k)=h,
$$
with respect to the variables $m$ and $h_1\ldots h_k$, $|h_i|<y$, $m\in I_k$, we find
\begin{equation}\label{summa-|Tn(alpha;x,y)|^2<sum r(h)e(alphah)}
|T(\alpha;x,y)|^{2^k}\leq(2y)^{2^k-k-1}\sum_{h}r(h)e(\alpha h),
\end{equation}
Note that if $h\neq 0$, then $r(h)\ll\tau_{k+1}(h)\ll h^{\varepsilon}$. From the inequality (\ref{formula gk(u;h1,...hk)>0}), it follows that the equation
$$
h_1\ldots h_kg_k(m;h_1,\ldots,h_k)=0
$$
has only solutions of the form $(0,h_2,\ldots,h_k,m)$, $(h_1,0,h_3,\ldots,h_k,m)$, $\ldots$, $(h_1,\ldots,h_{k-1},0, m)$, and the number of such solutions is estimated as
$$
r(0)\le k\sum_{|h_2|<y}\ldots\sum_{|h_k|<y}\sum_{m\in I_3}1\le k(2y)^{k-1}|I_k|\le k2^{k-1}y^k.
$$
On the other hand,
\begin{align}
\label{formula |Tn(alpha;x,y)|2k=sumrho(h)e(-alphah)}
|T(\alpha;x,y)|^{2^k}=\sum_{h}\rho(h)(-\alpha h),
\end{align}
where $\rho(h)$ is the number of solutions to the equation
$$
s_1^n+\ldots +s_\nu^n-t_1^n-\ldots -t_\nu^n=h,\qquad x-y<s_1,t_1,\ldots ,s_\nu,t_\nu\le x,\qquad \nu=2^k.
$$
In equality (\ref{formula |Tn(alpha;x,y)|2k=sumrho(h)e(-alphah)}), setting $\alpha=0$, we find
\begin{equation}\label{O srednem znach summ Weyl17}
\sum_{h}\rho(h)=|T(0;x,y)|^{2^k}\le y^{2^k}.
\end{equation}
Using the induction assumption, i.e., the relation (\ref{formula int|Tn(alpha;x,y)|2kdalpha<=y2k-k+varepsilon}), we have
$$
\rho(0)=\int_0^1 |T(\alpha;x,y)|^{2^k}d\alpha \le y^{2k-k+\varepsilon}.
$$
Multiplying (\ref{summa-|Tn(alpha;x,y)|^2<sum r(h)e(alphah)}) and (\ref{formula |Tn(alpha;x,y)|2k=sumrho(h)e(-alphah)}), integrating with respect to $\alpha$, and then using the values of $r(0)$, $\rho(0)$, the estimate $r(h)\ll h^\varepsilon$, and relation (\ref{O srednem znach summ Weyl17}), we find
\begin{align*}
\int_0^1|&T(\alpha;x,y)|^{2^{k+1}}d\alpha\le(2y)^{2^k-k-1}\int_0^1\sum_{h}r(h)e(\alpha h)\sum_{h'}\rho(h')e(-\alpha h')d\alpha=\\
&=(2y)^{2^k-k-1}\left(r(0)\rho(0)+\sum_{h\ne 0}r(h)\rho(h)\right)\le(2y)^{2^k-k-1}\left(r(0)\rho(0)+\max\limits_{h\ne 0}r(h)\sum_{h\ne 0}\rho(h)\right)\ll \\
&\ll y^{2^k-k-1}\left(y^k\cdot y^{2k-k+\varepsilon}+y^\varepsilon \cdot y^{2^k}\right)\ll y^{2^{k+1}-k-1+\varepsilon}.
\end{align*}

\section{Proof of Theorem \ref{TeorAsForEstKubPr}}
Without loss of generality, let's assume that
\begin{equation}\label{formula mu1<=...mu9}
H=N^{1-\theta(n,r)+\varepsilon},\qquad \theta(n,r)=\frac2{(n^2-n)(r+1)},\qquad \mu_1\le\ldots\le\mu_r.
\end{equation}
Using the notations
\begin{align*}
&N_k=(\mu_kN+H)^\frac1n, \qquad H_k=(\mu_kN+H)^\frac1n-(\mu_kN-H)^\frac1n,\qquad\tau=2(n-1)nN_1^{n-2}H_1,
\end{align*}
we express the number of solutions of the Diophantine equation (\ref{Formula-x1n+...+xrn=N}) under the conditions (\ref{formula |xin-N/n|<=H}) as
\begin{align*}
J_{n,r}(N,H)&=\int_{-\frac1\tau}^{1-\frac1\tau}e(-\alpha N)\prod_{k=1}^r\sum_{|m^n-\mu_kN|\le H}e(\alpha m^n)d\alpha= \int_{-\frac1\tau}^{1-\frac1\tau}e(-\alpha N)\prod_{k=1}^r\left(T(\alpha;N_k,H_k)+\theta_k\right)d\alpha,
\end{align*}
where $|\theta_k|$ is equal to $1$ if $N_k-H_k$ is an integer and $0$ otherwise. The upper bound of $N_k$ and the length of $H_k$ in the sum $T(\alpha;N_k,H_k)$ with respect to the parameters $N$ and $H$ are expressed through the following asymptotic formulas
\begin{align}
N_k&=\mu_k^\frac1nN^\frac1n\left(1+\frac{H}{\mu_kN}\right)^\frac1n=\mu_k^\frac1nN^\frac1n\left(1+O\left(\frac HN\right)\right), \label{formula poryadok Nk}\\
H_k&=\mu_k^\frac1nN^\frac1n\left(\left(1+\frac{H}{\mu_kN}\right)^\frac1n-\left(1-\frac{H}{\mu_kN}\right)^\frac1n\right)=
\frac{2H}{n\mu_k^{1-\frac1n}N^{1-\frac1n}}\left(1+O\left(\frac{H^2}{N^2}\right)\right). \label{formula poryadok Hk}
\end{align}
For $\nu=1,2,\ldots,r$ and $1\le i_1<\ldots<i_\nu\le r$, introducing the notation
$$
\mathscr{D}_\nu=\mathscr{D}(i_1,\ldots,i_\nu)=\{1,2,\ldots,r\}\setminus\{i_1,\ldots,i_\nu\}
$$
and using the identity
\begin{align*}
\prod_{k=1}^r&\left(T(\alpha;N_k,H_k)+\theta_k\right)=\prod_{k=1}^rT(\alpha;N_k,H_k)+\\
&+\sum_{\nu=1}^{r-1}\sum_{1\le i_1<\ldots<i_\nu\le r}\prod_{j=1}^\nu T(\alpha;N_{i_j},H_{i_j})
\prod_{k\in\mathscr{D}_\nu} \theta_k+\prod_{k=1}^r\theta_k,
\end{align*}
we represent $J_{n,r}(N,H)$ as
\begin{align}
&J_{n,r}(N,H)=\int_{-\frac1\tau}^{1-\frac1\tau}e(-\alpha N)\prod_{k=1}^rT(\alpha;N_k,H_k)d\alpha+R_{n,r}(N,H),\label{formula J_{n,r}(N,H)=int...+R1(n,H)}\\
&R_{n,r}(N,H)=\sum_{\nu=1}^{r-1}\sum_{1\le i_1<\ldots<i_\nu\le r}\prod_{k\in\mathscr{D}_\nu} \theta_k\int_{-\frac1\tau}^{1-\frac1\tau}e(-\alpha N)\prod_{j=1}^\nu T(\alpha;N_{i_j},H_{i_j})d\alpha. \nonumber
\end{align}
In the sum $R_1(N,H)$, going to estimates and using the fact that the geometric mean of non-negative numbers does not exceed their arithmetic mean, and then associating each $\nu$ with the number $k$ uniquely determined by the relation $2^k\le\nu<2^{k+1}$, remembering that $r-1=2^n$, we have
\begin{align}
R_{n,r}(N,H)&\le\sum_{\nu=1}^{r-1}\sum_{1\le i_1<\ldots<i_\nu\le r}\int_0^1\prod_{j=1}^\nu|T(\alpha;N_{i_j},H_{i_j})|d\alpha\le \nonumber\\
&\le\sum_{\nu=1}^{r-1}\frac{I(\nu)}\nu=I(1)+\sum_{k=1}^{n-1}\sum_{\nu=2^k}^{2^{k+1}-1}\frac{I(\nu)}\nu+\frac{I(r-1)}{r-1},\label{formula R1(N,H)} \\
&I(\nu)=\sum_{1\le i_1<\ldots<i_{\nu}\le r}\sum_{j=1}^\nu\int_0^1 \left|T(\alpha;N_{i_j},H_{i_j})\right|^\nu d\alpha.\nonumber
\end{align}
For $2\le\nu\le r-1$, let's estimate $I(\nu)$ from above. Using the trivial estimate of the sum $|T(\alpha;N_{i_j},H_{i_j})|$ and Theorem \ref{Teorema-obobsh-Hua}, we have
\begin{align*}
I(\nu)&\le\sum_{1\le i_1<\ldots<i_{\nu}\le r}\sum_{j=1}^\nu\max\left|T(\alpha;N_{i_j},H_{i_j})\right|^{\nu-2^k} \int_0^1\left|T(\alpha;N_{i_j},H_{i_j})\right|^{2^k} d\alpha\le\\
&\le\sum_{1\le i_1<\ldots<i_{\nu}\le r}\sum_{j=1}^\nu|H_{i_j}|^{\nu-2^k}\cdot|H_{i_j})|^{2^k-k+\varepsilon} \le\sum_{1\le i_1<\ldots<i_{\nu}\le r}\sum_{j=1}^\nu|H_{i_j}|^{\nu-k+\varepsilon}.
\end{align*}
Substituting this estimate into the right-hand side of (\ref{formula R1(N,H)}), and using the trivial estimate $J(1)$, and then utilizing the inequality $H_k\ll HN^{\frac1n-1}$, which follows from (\ref{formula poryadok Hk}), we find
\begin{align*}
R_{n,r}(N,H)&\ll \sum_{1\le i_1\le r}H_{i_1}+\sum_{k=1}^{n-1}\sum_{\nu=2^k}^{2^{k+1}-1}\sum_{1\le i_1<\ldots<i_{\nu}\le r}\sum_{j=1}^\nu|H_{i_j}|^{\nu-k+\varepsilon}+\sum_{1\le i_1<\ldots<i_{r-1}\le r}\sum_{j=1}^{r-1}|H_{i_j}|^{r-1-n+\varepsilon}\\
&\ll rHN^{\frac1n-1}+\sum_{k=1}^{n-1}\sum_{\nu=2^k}^{2^{k+1}-1}C_r^\nu\nu\left(HN^{\frac1n-1}\right)^{\nu-k+\varepsilon}+ C_r^{r-1}(r-1)\left(HN^{\frac1n-1}\right)^{r-1-n+\varepsilon}\ll\\
&\ll HN^{\frac1n-1}+\sum_{k=1}^{n-1}\left(HN^{\frac1n-1}\right)^{2^{k+1}-1-k+\varepsilon}+ \left(HN^{\frac1n-1}\right)^{r-1-n+\varepsilon}\ll\\
&\ll \left(\frac{H}{N^{1-\frac1n}}\right)^{r-1-n+\varepsilon}=\frac{H^{r-1}}{N^{r-\frac rn}} \left(\frac{N^{1-\frac1{n^2}+\frac{n-1}{n^2(n-\varepsilon)}\varepsilon}}H\right)^{n-\varepsilon}.
\end{align*}
Hence, considering that
$$
1-\frac1{n^2}+\frac{n-1}{n^2(n-\varepsilon)}\varepsilon\le1-\theta(n)+0,5\varepsilon,
$$
we find
\begin{align*}
R_{n,r}(N,H))&\ll\frac{H^{r-1}\lnc^\frac1n}{N^{r-\frac rn}} \left(\frac{N^{1-\theta(n)+\varepsilon}}H\right)^{n-\varepsilon}N^{-0,5\varepsilon(n-\varepsilon)} \ll\frac{H^{r-1}}{N^{r-\frac rn}\lnc}.
\end{align*}
From here and (\ref{formula J_{n,r}(N,H)=int...+R1(n,H)}), we have
\begin{align}
J_{n,r}(N,H)&=\int_{-\frac1\tau}^{1-\frac1\tau}e(-\alpha N)\prod_{k=1}^rT(\alpha;N_k,H_k)d\alpha +O\left(\frac{H^{r-1}}{N^{r-\frac rn}\lnc}\right).\label{formula J_{n,r}(N,H)=int...+O(H(r-r/n+varisilon))}
\end{align}
According to Dirichlet's theorem on the approximation of real numbers by rational numbers, every $\alpha$ from the interval $[-\tau^{-1}, 1-\tau^{-1}]$  can be represented as
\begin{equation}\label{formula alpha=a/q+lambda}
 \alpha =\frac{a}{q}+\lambda ,\qquad (a,q)=1,\qquad 1\le q\le \tau , \qquad |\lambda |\le \frac{1}{q\tau }.
\end{equation}
It is easy to see that in this representation $0\le a\le q-1$, and $a=0$ only when $q=1$. Let $\mathfrak{M}$ denote those $\alpha$ for which the condition $q\le H_r\lnc^{-1}$ holds in the representation (\ref{formula alpha=a/q+lambda}). Let $\mathfrak{m}$ denote the remaining $\alpha$. The set $\mathfrak{M}$ consists of non-overlapping intervals. We divide the set $\mathfrak{M}$ into the sets $\mathfrak{M}_1$ and $\mathfrak{M}_2$:
\begin{align}
&\mathfrak{M}_1=\bigcup_{1\le q\le H_r\lnc^{-1}}\bigcup_{\substack{a=0\\ (a,q)=1}}^{q-1}\mathfrak{M}_1(a,q),\qquad
\mathfrak{M}_1(a,q)=\left[\frac{a}{q}-\eta_q \le\alpha\le\frac{a}{q}+\eta_q \right],\nonumber\\ &\mathfrak{M}_2=\mathfrak{M}\setminus\mathfrak{M}_1,\qquad \eta_q=\frac1{2nqN_r^{n-1}},\qquad \eta=\frac{\lnc}{2nH_rN_{r-1}}\le\eta_q.\label{formula opr eta i etaq}
\end{align}
Denoting by $J(\mathfrak{M}_1)$, $J(\mathfrak{M}2)$, and $J(\mathfrak{m})$ the integrals over the sets $\mathfrak{M}1$, $\mathfrak{M}2$, and $\mathfrak{m}$, respectively, taking into account (\ref{formula J_{n,r}(N,H)=int...+O(H(r-r/n+varisilon))}), we obtain
\begin{equation}\label{formula J_{3,9}(N,H)=I(M1)+I(M2)+I(m)+O(..)}
J_{n,r}(N,H)=J(\mathfrak{M}_1)+J(\mathfrak{M}_2)+J(\mathfrak{m})+O\left(\frac{H^{r-1}}{N^{r-\frac rn}\lnc}\right).
\end{equation}
In the last formula, the first term, i.e., $J(\mathfrak{M}1)$, contributes the main term to the asymptotic formula for $J_{n,r}(N,H)$, while $J(\mathfrak{M}_2)$ and $J(\mathfrak{m})$ are part of its residual term.

\subsection{Calculation of the integral $J(\mathfrak{M}_1)$}
By the definition of the integral $J(\mathfrak{M}_1)$, we have:
\begin{align}\label{formula I(M1)=1}
&J(\mathfrak{M}_1)=\hspace{-14pt}\sum_{q\le H_r\lnc^{-1}}\sum_{\substack{a=0\\(a,q)=1}}^{q-1}\ \int\limits_{|\lambda|\le\eta_q} \prod_{k=1}^rT\left(\frac{a}{q}+\lambda;N_k,H_k\right)e\left(-\left(\frac{a}{q}+\lambda\right)N\right)d\lambda.
\end{align}
%ƒл€ суммы $T\left(\frac{a}{q}+\lambda;N_k,H_k\right)$ выполн€ютс€ оба услови€ следстви€ \ref{Sledst1 Teor ob poved kor trig summi Weyl} теоремы \ref{Teor ob poved kor trig summi Weyl}. ƒействительно, ввиду соотношений (\ref{formula poryadok Nk}), (\ref{formula poryadok Hk}) и (\ref{formula mu1<=...mu9}) имеет место неравенство
For the sum $T\left(\frac{a}{q}+\lambda;N_k,H_k\right)$ both conditions of the corollary \ref{Sledst1 Teor ob poved kor trig summi Weyl} of theorem \ref{Teor ob poved kor trig summi Weyl} are satisfied. Indeed, in the view of the relations (\ref{formula poryadok Nk}), (\ref{formula poryadok Hk}) and (\ref{formula mu1<=...mu9}) the following inequality holds
\begin{equation}\label{formula tau=12N1H1>=12NrHr}
\begin{split}
\tau=&2(n-1)nN_1^{n-2}H_1=\frac{4(n-1)H}{\mu_1^\frac1nN^\frac1n}\left(1+O\left(\frac HN\right)\right)\ge\\
 &\ge\frac{4(n-1)H}{\mu_k^\frac1nN^\frac1n}\left(1+O\left(\frac HN\right)\right)=2(n-1)nN_k^{n-2}H_k,
\end{split}
\end{equation}
and from the relations $|\lambda|\le\eta_q$, $\eta_q=\dfrac1{2nqN_r^{n-1}}$, and $\dfrac1{2nqN_r^{n-1}}\le\dfrac1{2qN_k^{n-1}}$, it follows that
$$
|\lambda|\le\frac1{2nqN_r^{n-1}}.
$$
So, according to this corollary, for $k=1,\ldots,r$, we have
\begin{align*}
T\left(\frac aq+\lambda,N_k,H_k\right)=\frac{H_kS(a,q)}{q}\gamma(\lambda;N_k,H_k)+R,\qquad R\ll q^{\frac12+\varepsilon}.
\end{align*}
Multiplying both sides of these formulas for all $k=1,2,\ldots,r$ and then applying the identity
\begin{align*}
&\prod_{k=1}^r(a_k+b)=\prod_{k=1}^ra_k+b^r+\sum_{\nu=1}^{r-1}b^{r-\nu}\sum_{1\le i_1<\ldots<i_\nu\le r}\prod_{k=1}^\nu a_{i_k},
\end{align*}
with $a_k=\dfrac{H_kS(a,q)}{q}\gamma(\lambda;N_k,H_k)$ and $b=R$, we obtain
\begin{align*}
\prod_{k=1}^rT\left(\frac aq+\lambda,N_k,H_k\right)&=\frac{S^r(a,q)}{q^r}\prod_{k=1}^rH_k\gamma(\lambda;N_k,H_k)+R^r+\\
&+\sum_{\nu=1}^{r-1}R^{r-\nu}\sum_{1\le i_1<\ldots<i_\nu\le r}\prod_{k=1}^\nu\frac{H_{i_k}S(a,q)}q\gamma(\lambda;N_{i_k},H_{i_k}).
\end{align*}
Using the relations $R\ll q^{\frac12+\varepsilon}$ and $|S(a,q)|\ll q^\frac{n-1}n$ (lemma \ref{Lemma otsenka polnoy trigsumm}), we estimate the last two terms from above:
\begin{align*}
\prod_{k=1}^rT&\left(\frac{a}{q}+\lambda,N_k,H_k\right)-\frac{S^r(a,q)}{q^r}\prod_{k=1}^rH_k\gamma(\lambda;N_k,H_k)\ll\\
&\ll\sum_{\nu=1}^{r-1}q^{\left(\frac12+\varepsilon\right)(r-\nu)-\frac{\nu}n}\sum_{1\le i_1<\ldots<i_\nu\le r}\prod_{k=1}^\nu H_{i_k} |\gamma(\lambda;N_{i_k},H_{i_k})|+q^{0,5r+r\varepsilon}.
\end{align*}
From here and from the formula (\ref{formula I(M1)=1}), we find
\begin{align}
J(\mathfrak{M}_1)&=\prod_{i=1}^rH_i\hspace{-10pt}\sum_{q\le H_r\lnc^{-1}}\hspace{-12pt}\mathscr{A}(r,q)\hspace{-7pt}\sum_{\substack{a=0\\(a,q)=1}}^{q-1}\hspace{-8pt}\frac{S^r(a,q)}{q^r} e\left(-\frac{aN}q\right)+R_1(\mathfrak{M}_1)+R_2(\mathfrak{M}_1), \label{formula I(M1)-1}\\
\mathscr{A}(r,q)&=\int_{|\lambda|\le\eta_q }\prod_{k=1}^r\gamma(\lambda;N_k,H_k)e(-\lambda N)d\lambda,\nonumber\\
R_1(\mathfrak{M}_1)&\ll \sum_{\nu=1}^{r-1}\sum_{q\le H_r\lnc^{-1}}q^{\sigma(\nu)}\hspace{-8pt}\sum_{1\le i_1<\ldots<i_\nu\le r} \int\limits_{|\lambda|\le\eta_q}\prod_{k=1}^\nu H_{i_k} |\gamma(\lambda;N_{i_k},H_{i_k})|d\lambda,\label{formula R1(M1)-1}\\
\sigma(\nu)&= \left(\frac12+\varepsilon\right)r+1-\left(\frac12+\frac1n+\varepsilon\right)\nu,\nonumber\\
R_2(\mathfrak{M}_1)&\ll\sum_{q\le H_r\lnc^{-1}}q^{0,5r+r\varepsilon}\varphi(q)\cdot2\eta_q =\frac{1}{nN_r^{n-1}}\sum_{q\le H_r\lnc^{-1}}\varphi(q)q^{0,5r-1+r\varepsilon}.\nonumber
\end{align}

\subsection{Estimation of $R_2(\mathfrak{M}_1)$}
Using formulas (\ref{formula poryadok Hk}) and (\ref{formula poryadok Nk}), we have
\begin{align*}
R_2(\mathfrak{M}_1)&\ll\frac1{N_r^{n-1}}\left(\frac{H_r}{\lnc}\right)^{0,5r+1+r\varepsilon}\ll N^\frac{1-n}n\left(\frac{H}{N^\frac{n-1}n\lnc}\right)^{0,5r+1+r\varepsilon}=\\
&=\frac{H^{r-1}}{N^{r-\frac rn}\lnc^{0,5r+1+r\varepsilon}} \left(\frac{N^{1-\frac1n}}H\right)^{\frac r2-2-r\varepsilon}.
\end{align*}
%ќтсюда име€ в виду, что
Hence, considering that
$$
1-\frac1n<1-\theta(n),\qquad \frac r2-2-r\varepsilon>2,
$$
and also using the relation $H=N^{1-\theta(n,r)+\varepsilon}$, we obtain
\begin{align}
R_2(\mathfrak{M}_1)&\ll\frac{H^{r-1}}{N^{r-\frac rn}\lnc^{0,5r+1+r\varepsilon}}\left(\frac{N^{1-\theta(n)+\varepsilon}}H\right)^{\frac r2-2-r\varepsilon}N^{-\varepsilon\left(\frac r2-2-r\varepsilon\right)}\ll\frac{H^{r-1}}{N^{r-\frac rn}\lnc^{r-1}}.
\label{formula R2(M1)}
\end{align}

\subsection{Estimation of $R_1(\mathfrak{M}_1)$} First, let's estimate the exponential integral
\begin{align*}
&\gamma(\lambda;N_{i_k},H_{i_k})=\int_{-0,5}^{0,5}e\left(f_{i_k}(u)\right)du, \qquad f_{i_k}(u)=\lambda\left(N_{i_k}-\frac{H_{i_k}}2+H_{i_k}u\right)^n.
\end{align*}
Utilizing the relations (\ref{formula poryadok Nk}) and (\ref{formula poryadok Hk}), for $i=1,\ldots,n$, we obtain
\begin{equation}\label{formula poryadok Nk(n-i)Hki i Hk/Nk}
\begin{split}
&N_k^{n-i}H_k^i=\frac{2^iN^{1-i}H^i}{n\mu_k^{i-1}}\left(1+O\left(\frac HN\right)\right), \qquad\frac{H_k}{N_k}=\frac{2H}{n\mu_kN}\left(1+O\left(\frac HN\right)\right).
 \end{split}
\end{equation}
Using these relations, let's estimate the lower bound for $f_{i_k}'(u)$. We have
\begin{align*}
|f_{i_k}'(u))|&=n|\lambda|H_{i_k}\left(N_{i_k}-\frac{H_{i_k}}2+H_{i_k}u\right)^{n-1}
\ge n|\lambda|H_{i_k}(N_{i_k}-H_{i_k})^{n-1}=\\
&=n|\lambda|N_{i_k}^{n-1}H_{i_k}\left(1-\frac{H_{i_k}}{N_{i_k}}\right)^{n-1}=2|\lambda|H\left(1+O\left(\frac HN\right)\right)\ge |\lambda|H.
\end{align*}
From here and using lemma \ref{Lemma otsenka trig integ po pervoy proizv}, we obtain
\begin{equation}\label{formula otsenka int gamma3(...)}
|\gamma(\lambda;N_{i_k},H_{i_k})|\le\min\left(1,\frac1{|\lambda|H}\right).
\end{equation}
Substituting this estimate into the right-hand side of (\ref{formula R1(M1)-1}), we obtain
\begin{align}
R_1(\mathfrak{M}_1)\ll&\sum_{\nu=1}^{r-1}\sum_{q\le H_r\lnc^{-1}}q^{\sigma(\nu)}\sum_{1\le i_1<\ldots<i_\nu\le r} \int_{|\lambda|\le\eta_q}\prod_{k=1}^\nu H_{i_k}\min\left(1,\frac1{|\lambda|H}\right)d\lambda=\nonumber\\
&=2\sum_{\nu=1}^{r-1}I(\nu)\sum_{q\le H_r\lnc^{-1}}q^{\sigma(\nu)}\sum_{1\le i_1<\ldots<i_\nu\le r}\prod_{k=1}^\nu H_{i_k},\label{formula R1(M1)-2}\\
&\qquad I(\nu)=\int_0^{\eta_q}\min\left(1,\frac1{\lambda^\nu H^\nu}\right)d\lambda.\nonumber
\end{align}
Using of the condition $q \leq H_r\lnc^{-1}$, and then using the relation (\ref{formula poryadok Nk(n-i)Hki i Hk/Nk}), we find
\begin{align*}
\eta_qH=\frac{H}{2nqN_r^{n-1}}\ge\frac{H\lnc}{2nH_rN_r^{n-1}}=\frac{H\lnc}{2n\cdot \frac{2H}n\left(1+O\left(\frac HN\right)\right)} =\frac{\lnc}{4\left(1+O\left(\frac HN\right)\right)}\ge\frac{\lnc}5,
\end{align*}
that is, $H^{-1} < \eta_q$. For $\nu \geq 2$, by dividing the integration interval in the integral $I(\nu)$ into the intervals $[0, H^{-1}]$ and $[H^{-1}, \eta_q]$, we have
\begin{align*}
I(\nu)&=\int_0^{H^{-1}}d\lambda+\frac1{H^\nu}\int_{H^{-1}}^{\eta_q }\frac{d\lambda}{\lambda^\nu}
=\frac1H\left(1+ \frac1{\nu-1}\left(1-\left(\frac1{\eta_qH}\right)^{\nu-1}\right)\right)\le\frac{\nu}{(\nu-1)H}.
\end{align*}
In the case of $\nu=1$, similarly, we obtain
\begin{align*}
I(1)&=\int_0^{H^{-1}}d\lambda+\frac1{H}\int_{H^{-1}}^{\eta_q }\frac{d\lambda}{\lambda}=\frac1H+ \frac{\ln(\eta_qH)}H\le\frac{\lnc}{H}.
\end{align*}
Using the relationships (\ref{formula poryadok Hk}) and (\ref{formula mu1<=...mu9}), we find
\begin{align*}
&\prod_{k=1}^\nu H_{i_k}=\prod_{k=1}^\nu \frac{2H}{n\mu_{i_k}^{1-\frac1n}N^{1-\frac1n}}\left(1+O\left(\frac{H^2}{N^2}\right)\right)\ll \left(\frac{H}{N^{1-\frac1n}}\right)^\nu.
\end{align*}
Substituting the right-hand side of this inequality and the estimate for the integral $I(\nu)$ into (\ref{formula R1(M1)-2}), and then using the formula (\ref{formula poryadok Hk}) for $k=r$, we have:
\begin{align*}
R_1(\mathfrak{M}_1)&\ll\sum_{\nu=1}^{r-1}\sum_{1\le i_1<\ldots<i_\nu\le r}\frac{\lnc^{1+\nu}}H\left(\frac{H}{N^{1-\frac1n}\lnc}\right)^\nu\sum_{q\le H_r\lnc^{-1}}q^{\sigma(\nu)}\ll\\
&\ll\sum_{\nu=1}^{r-1}C_r^\nu\frac{\lnc^{1+\nu}}H\left(\frac{H}{N^{1-\frac1n}\lnc}\right)^{\sigma(\nu)+\nu+1}\ll \frac{\lnc^r}H\left(\frac{H}{N^{1-\frac1n}\lnc}\right)^{\sigma(r-1)+r}=\\
&=\frac{H^{r-1}}{N^{r-\frac rn}}\lnc^{\frac rn-\frac32-\frac1n-\varepsilon}\left(\frac{N^{1-\frac1n}}H\right)^{\frac rn-\frac32-\frac1n-\varepsilon}.
\end{align*}
Therefore, considering that
$$
1-\frac1n<1-\theta(n),\qquad \frac rn-\frac32-\frac1n-\varepsilon>1,
$$
and using the relation $H=N^{1-\theta(n,r)+\varepsilon}$, we find
\begin{align}
R_1(\mathfrak{M}_1)&\ll\frac{H^{r-1}}{N^{r-\frac rn}}\lnc^{\frac rn-\frac32-\frac1n-\varepsilon}\left(\frac{N^{1-\theta(n,r)+\varepsilon}} H\right)^{\frac rn-\frac32-\frac1n-\varepsilon} N^{-\varepsilon\left(\frac rn-\frac32-\frac1n-\varepsilon\right)}=\nonumber\\
&=\frac{H^{r-1}}{N^{r-\frac rn}}\lnc^{\frac rn-\frac32-\frac1n-\varepsilon}\cdot N^{-\varepsilon\left(\frac rn-\frac32-\frac1n-\varepsilon\right)}\ll\frac{H^{r-1}}{N^{r-\frac rn}\lnc^{r-1}}.\label{formula R1(M1)-3}
\end{align}

\subsection{Calculation of the integral $\mathscr{A}(r,q)$} We have
$$
\mathscr{A}(r,q)=\int_{|\lambda|\le\eta_q}\prod_{k=1}^r\gamma(\lambda;N_k,H_k)e(-\lambda N)d\lambda,\qquad \eta_q=\frac1{2nqN_r^{n-1}}.
$$
By dividing the integration interval into intervals $|\lambda|\le\eta$ and $\eta<|\lambda|\le\eta_q$, where the value $\eta$ is determined in formula (\ref{formula opr eta i etaq}), and denoting the integrals over these intervals as $\mathscr{A}_1(r,q)$ and $\mathscr{A}_2(r,q)$ respectively, we obtain
\begin{align}
&\mathscr{A}(r,q)=\mathscr{A}_1(r,q)+\mathscr{A}_2(r,q).\label{formula opr A(M1)-1}
\end{align}
Notice that the value of $\eta$ is determined in formula (\ref{formula opr eta i etaq}), and according to (\ref{formula poryadok Nk(n-i)Hki i Hk/Nk}), we have
\begin{equation}\label{formula poryadok eta}
\eta=\frac{\lnc}{2nH_rN_r^{n-1}}=\frac{\lnc}{4H\left(1+O\left(\frac HN\right)\right)}\le\frac{\lnc}H.
\end{equation}
In formula (\ref{formula opr A(M1)-1}), the integral $\mathscr{A}_1(r,q)$ provides the main term of the asymptotic formula for $\mathscr{A}(r,q)$, while $\mathscr{A}_2(r,q)$ contributes to the remainder term. Let's first find the asymptotic formula for $\mathscr{A}_1(r,q)$. Using the relation $N_k^n=\mu_kN+H$, formulas (\ref{formula poryadok Nk(n-i)Hki i Hk/Nk}), and (\ref{formula poryadok eta}), we have
\begin{align*}
&f_k(u)=\lambda\left(N_k+H_k\left(u-\frac12\right)\right)^n=\lambda N_k^n+\lambda\sum_{i=1}^nC_n^iN_k^{n-i}H_k^i\left(u-\frac12\right)^i=\\
&=\mu_kN\lambda+2H\lambda u+O\left(\frac{H^2}N|\lambda|\right)+ \lambda\sum_{i=2}^nC_n^i\frac{2^iN^{1-i}H^i}{n\mu_k^{i-1}}\left(u-\frac12\right)^i\left(1+O\left(\frac HN\right)\right)=\\
&=\mu_kN\lambda+2H\lambda u+R_3(N,H),\qquad R_3(N,H)\ll\frac{H^2}N\eta\ll\frac{H\lnc}N.
\end{align*}
From here, considering that $e(R_3(N,H))=1+O(HN^{-1}\lnc)$, we find
\begin{align*}
&e\left(f_k(u)\right)=e\left(\mu_kN\lambda\right)e(2Hu\lambda)+O\left(\frac{H\lnc}N\right).
\end{align*}
Therefore,
\begin{align*}
&\gamma(\lambda;N_k,H_k)=e\left(\mu_kN\lambda\right)\frac{\sin\left(2\pi H\lambda\right)}{2\pi H\lambda}+R_4(N,H),&R_4(N,H)\ll\frac{H\lnc}N.
\end{align*}
Multiplying both sides of these formulas for all $k=1,2,\ldots,r$, applying the identity
\begin{align*}
&\prod_{k=1}^r(a_k+b)=\prod_{k=1}^ra_k+b^r+\sum_{\nu=1}^{r-1}b^{r-\nu}\sum_{1\le i_1<\ldots<i_\nu\le r}\prod_{k=1}^\nu a_{i_k},
\end{align*}
and then using the condition $\mu_1+\mu_2+\dots+\mu_r=1$, we have:
\begin{align*}
\prod_{k=1}^r\gamma(\lambda;N_k,H_k)&=\prod_{k=1}^r\left(e\left(\mu_kN\lambda\right)\frac{\sin\left(2\pi H\lambda\right)}{2\pi H\lambda}+R_4(N,H)\right)=\\
&=\frac{\sin^r\left(2\pi H\lambda\right)}{(2\pi H\lambda)^r}e\left(N\lambda\right)+R_5(N,H),
\end{align*}
where
\begin{align*}
R_5(N,&H)=R_4^r(N,H) +\sum_{\nu=1}^{r-1}R_4^{r-\nu}(N,H)\sum_{1\le i_1<\ldots<i_\nu\le r}\prod_{k=1}^\nu e\left(\mu_{i_k}N\lambda\right)\frac{\sin\left(2\pi H\lambda\right)}{2\pi H\lambda}=\\
&=R_4^r(N,H)+\sum_{\nu=1}^{r-1}\frac{\sin^\nu\left(2\pi H\lambda\right)}{(2\pi H\lambda)^\nu}R_4^{r-\nu}(N,H)\sum_{1\le i_1<\ldots<i_\nu\le r}e\left(\lambda N\sum_{k=1}^\nu\mu_{i_k}\right)\ll\\
&\ll\sum_{\nu=0}^{r-1}\frac{|\sin(2\pi H\lambda)|^\nu}{|2\pi\lambda H|^\nu}\left(\frac{H\lnc}N\right)^{r-\nu} \ll
\frac{H^r\lnc^r}{N^r}+\frac{|\sin(2\pi H\lambda)|^{r-1}}{|2\pi\lambda H|^{r-1}}\cdot\frac{H\lnc}N.
\end{align*}
From here, using the definition of the integral $\mathscr{A}_1(k,q)$ and (\ref{formula poryadok eta}), we find:
\begin{align*}
\mathscr{A}_1(r,q)&=\int_{|\lambda|\le\eta}\left(\frac{\sin^r\left(2\pi H\lambda\right)}{(2\pi H\lambda)^r}e\left(N\lambda\right)+ R_5(N,H)\right) e(-\lambda N)d\lambda=\\
&=\frac1{\pi H}\int_0^{2\pi H\eta}\frac{\sin^rt}{t^r}dt+R_6(N,H)+R_7(N,H),\\
R_6(N,H)&\ll\int_{|\lambda|\le\eta}\frac{H^r\lnc^r}{N^r}d\lambda\le\frac{2H^{r-1}\lnc^{r+1}}{N^r}=\frac1{H\lnc^7}\cdot \frac{H^r\lnc^{r+8}}{N^r} \ll\frac1{H\lnc^7},\\
R_7(N,H)&\ll\frac{H\lnc}N\int_{|\lambda|\le\eta}\frac{\sin^{r-1}\left(2\pi H\lambda\right)}{(2\pi H\lambda)^{r-1}}d\lambda=
\frac{\lnc}{\pi N}\int_0^{2\pi H\eta}\frac{\sin^{r-1}t}{t^{r-1}}dt\ll\frac1{H\lnc^7}.
\end{align*}
By replacing the integral over $t$ with a nearby improper integral independent of $2\pi H\eta$ and utilizing the formula (\ref{formula poryadok Nk(n-i)Hki i Hk/Nk}) for $i=1$, we obtain:
\begin{align*}
&\mathscr{A}_1(r,q)=\frac1{\pi H}\int_0^\infty\frac{\sin^rt}{t^r}dt+R_8(N,H)+O\left(\frac1{H\lnc^7}\right),\\
&R_8(N,H)=\frac1{\pi H}\int_{2\pi H\eta}^\infty\frac{\sin^rt}{t^r}dt\le\frac1{\pi H}\cdot\frac1{(2\pi H\eta)^r}\ll \frac1{H^{r+1}\eta^r}=\\
&=\frac1{H^{r+1}}\cdot\left(\frac{2nN_r^{n-1}H_r}\lnc\right)^r\ll\frac1{H^{r+1}\lnc^r}\left(\frac{2H}n\right)^r\ll\frac1{H\lnc^r}.
\end{align*}
Using the lemma \ref{Lemma tochnoe znach trig integ} to calculate the improper integral with $m=1$, we have:
\begin{align}
&\mathscr{A}_1(r,q)=\frac{\gamma(n,r)}{H}+O\left(\frac1{H\lnc^7}\right).\label{formula opr A1(M1)-3}
\end{align}
\begin{align*}
\gamma(n,r)=&\frac{r^{r-1}-\frac{r}{1!}(r-2)^{r-1}+\frac{r(r-1)}{2!}(r-4)^{r-1}-\frac{r(r-1)(r-2)}{3!}(r-6)^{r-1}+\ldots}{2^r(r-1)!}.
\end{align*}
Now let's estimate the integral $\mathscr{A}_2(r,q)$ from above. Using the estimate from (\ref{formula otsenka int gamma3(...)}), and then applying the relation (\ref{formula poryadok eta}), we obtain:
\begin{align*}
\mathscr{A}_2(r,q)&\le 2\int_\eta^{\eta_q}\prod_{k=1}^r|\gamma(\lambda;N_k,H_k)|d\lambda\le2\int_\eta^{\eta_q}\min\left(1,\frac1{\lambda^rH^r}\right)d\lambda=\\
&=\frac2{H^r}\int_\eta^{\eta_q}\frac{d\lambda}{\lambda^r}=\frac2{(r-1)H^r}\left(\frac1{\eta^{r-1}}-\frac1{\eta_q^{r-1}}\right)\le
\frac2{(r-1)H^r\eta^{r-1}}=\\
&=\frac2{(r-1)H^r}\cdot\left(\frac{4H\left(1+O\left(\frac HN\right)\right)}{\lnc}\right)^{r-1}\ll\frac1{H\lnc^{r-1}}.
\end{align*}
From this estimate and the formula (\ref{formula opr A1(M1)-3}), considering (\ref{formula opr A(M1)-1}), we find:
\begin{align}
&\mathscr{A}(r,q)=\frac{\gamma(n,r)}{H}+O\left(\frac1{H\lnc^{r-1}}\right).\label{formula opr A(M1)-2}
\end{align}

\subsection{Derivation of the Asymptotic Formula for the Integral $J(\mathfrak{M}_1)$}
Substituting the right-hand sides of the formulas (\ref{formula opr A(M1)-2}), (\ref{formula R1(M1)-3}), and (\ref{formula R2(M1)}) into (\ref{formula I(M1)-1}), we obtain:
\begin{align}
J(\mathfrak{M}_1)&=\frac{\gamma(n,r)}{H}\mathfrak{S}\left(N,\frac{H_r}{\lnc}\right)\prod_{i=1}^rH_i +R_9(N,H)+O\left(\frac{H^{r-1}}{N^{r-\frac rn}\lnc^{r-1}}\right), \label{formula I(M1)-2}\\
&\mathfrak{S}\left(N,\frac{H_r}{\lnc}\right)=\sum_{q\le H_r\lnc^{-1}}\sum_{\substack{a=0\\(a,q)=1}}^{q-1}\frac{S^r(a,q)}{q^r} e\left(-\frac{aN}q\right), \nonumber\\
&R_9(N,H)\ll\frac1{H\lnc^{r-1}}\prod_{i=1}^rH_i\sum_{q\le H_r\lnc^{-1}}\left|\sum_{\substack{a=0\\(a,q)=1}}^{q-1}\frac{S^r(a,q)}{q^r} e\left(-\frac{aN}q\right)\right|.\nonumber
\end{align}
Let's compute the double sum $\mathfrak{S}\left(N,H_r\lnc^{-1}\right)$. To do this, we will replace the sum over $q$ with a closely related infinite series $\mathfrak{S}(N)$, independent of $H_r\lnc^{-1}$. Utilizing Lemma \ref{Lemma otsenka polnoy trigsumm}, the relations (\ref{formula poryadok Nk}), and (\ref{formula poryadok Hk}), and then substituting explicit values for the parameters $H$ and $r$, namely using formula (\ref{formula mu1<=...mu9}), we have:
\begin{align*}
&\left|\sum_{q>H_r\lnc^{-1}}\sum_{\substack{a=0\\(a,q)=1}}^{q-1}\frac{S^r(a,q)}{q^r} e\left(-\frac{aN}q\right)\right|\ll \sum_{q>H_r\lnc^{-1}}q^{-\frac rn+1}\ll\\
&\hspace{70pt}\ll\left(\frac{H_r}{\lnc}\right)^{-\frac rn+2}\ll\left(\frac{N^{1-\frac1n}\lnc}H\right)^{\frac rn-2} \ll\frac1{\lnc^{r-1}}.
\end{align*}
Hence
\begin{align}\label{formula osobi ryad}
&\mathfrak{S}\left(N,\frac{H_r}{\lnc}\right)=\mathfrak{S}(N)+O\left(\frac1{\lnc^{r-1}}\right),\\ &\mathfrak{S}(N)=\sum_{q=1}^\infty\sum_{\substack{a=0\\(a,q)=1}}^{q-1}\frac{S^r(a,q)}{q^9} e\left(-\frac{aN}q\right).\nonumber
\end{align}
The sum of the special series $\mathfrak{S} (N)$ exceeds some positive number $c(N)$ (see \cite{Vaughan-1985}, theorem 4.6).

Applying the formula (\ref{formula poryadok Hk}), we have:
\begin{align}\label{formula H1...H9}
\prod_{i=1}^rH_i=\frac{2^rH^r}{n^rN^{r-\frac rn}}\prod_{i=1}^r\mu_i^{-1+\frac1n}\left(1+O\left(\frac{H^2}{N^2}\right)\right).
\end{align}
For the estimation of $R_9(N,H)$, using the last formula and Lemma \ref{Lemma otsenka polnoy trigsumm} along with the condition $\dfrac rn-1\ge2$, we obtain:
\begin{align}\label{formula R8(N,H)}
R_9(N,H)&\ll\frac{H^{r-1}}{N^{r-\frac rn}\lnc^{r-1}}\sum_{q\le H_r\lnc^{-1}}\frac1{q^{\frac rn-1}}\ll\frac{H^{r-1}}{N^{r-\frac rn}\lnc^{r-1}}.
\end{align}
Substituting (\ref{formula osobi ryad}), (\ref{formula H1...H9}), and (\ref{formula R8(N,H)}) into (\ref{formula I(M1)-2}), we find:
\begin{align}
J(\mathfrak{M}_1)&=\frac{2^r\gamma(n,r)}{n^r}\prod_{i=1}^r\mu_i^{-1+\frac1n}\mathfrak{S}(N)\frac{H^{r-1}}{N^{r-\frac rn}}+ O\left(\frac{H^{r-1}}{N^{r-\frac rn}\lnc^{r-1}}\right)
. \label{formula I(M1)-3}
\end{align}

\subsection{Estimation of the integral $J(\mathfrak{M}_2)$}
%»меем
We have
\begin{align}\label{formula I(M2)-1}
&J(\mathfrak{M}_2)=\int_{\mathfrak{M}_2}\prod_{k=1}^rT(\alpha;N_k,H_k)e(-\alpha N)d\alpha.
\end{align}
The sums $T(\alpha;N_k,H_k)$ in the product $\prod\limits_{k=1}^rT(\alpha;N_k,H_k)$ are symmetric. Without loss of generality, we can assume that the following relation holds:
$$
\max_{1\le k\le r}\max_{\alpha\in \mathfrak{M}_2}\left|T(\alpha;N_k,H_k)\right|=\max_{\alpha\in \mathfrak{M}_2}\left|T(\alpha;N_\nu,H_\nu)\right|,\qquad 1\le\nu\le r.
$$
With this equality in mind, transitioning to estimates in the integral (\ref{formula I(M2)-1}), and using the fact that the geometric mean of non-negative numbers does not exceed their arithmetic mean, and then applying Theorem \ref{Teorema-obobsh-Hua} and the relation $H_\nu\le H_1\ll HN^{-1+\frac1n}$, we have:
\begin{align}
J(\mathfrak{M}_2)&\le\max_{\alpha\in \mathfrak{M}_2}\left|T(\alpha;N_\nu,H_\nu)\right|\int_0^1\prod_{\substack{k=1\\k\neq \nu}}^r\left|T(\alpha;N_k,H_k)\right|d\alpha\le\nonumber\\
&\le\max_{\alpha\in \mathfrak{M}_2}\left|T(\alpha;N_\nu,H_\nu)\right|\frac1{r-1}\sum_{\substack{k=1\\k\neq \nu}}^r\int_0^1\left|T(\alpha;N_k,H_k)\right|^{r-1}d\alpha\ll\nonumber\\
&\ll\max_{\alpha\in \mathfrak{M}_2}\left|T(\alpha;N_\nu,H_\nu)\right|\frac1{r-1}\sum_{\substack{k=1\\k\neq \nu}}^r H_k^{2^n-n+\varepsilon}\ll\nonumber\\
&\ll\left(\frac{H}{N^{1-\frac1n}}\right)^{2^n-n+\varepsilon}\hspace{-10pt}\max_{\alpha\in \mathfrak{M}_2}\left|T(\alpha;N_\nu,H_\nu)\right|=\nonumber\\
&=\frac{H^{r-1}}{N^{r-\frac rn}}\cdot\frac{N^{n-\frac1n-\frac{n-1}n\varepsilon}}{H^{n-\varepsilon}}\max_{\alpha\in \mathfrak{M}_2}\left|T(\alpha;N_\nu,H_\nu)\right|. \label{formula I(M2)-2}
\end{align}

Let's estimate $T(\alpha;N_\nu,H_\nu)$ for $\alpha$ from the set $\mathfrak{M}_2$. If $\alpha\in\mathfrak{M}_2$, then:
$$
\alpha=\frac aq+\lambda,\quad(a,q)=1,\quad\eta_q<|\lambda |\le\frac1{q\tau},\quad 1\le q\le\frac{H_r}\lnc, \quad\eta_q=\frac1{2nqN_r^{n-1}}.
$$
Let's consider two possible cases: $\eta_q < |\lambda| \leq \frac{1}{2nqN_\nu^{n-1}}$ and $\frac{1}{2nqN_\nu^{n-1}} < |\lambda| \leq \frac{1}{q\tau}$.

{\bf Case 1.} For the sum $T(\alpha;N_\nu,H_\nu)$, according to the relation (\ref{formula tau=12N1H1>=12NrHr}), the next inequality holds:
$$
\tau=2(n-1)nN_1^{n-2}H_1\ge2(n-1)nN_\nu^{n-2}H_\nu,
$$
that is the first condition of Corollary \ref{Sledst1 Teor ob poved kor trig summi Weyl} of Theorem \ref{Teor ob poved kor trig summi Weyl}, and the second condition follows from the case
$$
|\lambda |\le\frac1{2nqN_r^{n-1}}.
$$
According to this corollary, we have
\begin{align}\label{formula T3(lambda;Nr,Hr)}
T(\alpha;N_\nu,H_\nu)=\frac{H_\nu S(a,q)}{q}\gamma(\lambda;N_\nu,H_\nu)+O\left(q^{\frac12+\varepsilon}\right).
\end{align}
Estimating the exponential integral $\gamma(\lambda;N_\nu,H_\nu)$ using the estimate (\ref{formula otsenka int gamma3(...)}), we find
$$
\gamma(\lambda;N_\nu,H_\nu)\le\min\left(1,\frac1{H|\lambda|}\right)\le\frac1{H\eta_q}= \frac{2nqN_r^{n-1}}H\ll\frac{qN^\frac{n-1}n}{H}.
$$
Substituting this estimate and the estimate for the sum $S(a,q)$ from Lemma \ref{Lemma otsenka polnoy trigsumm} into (\ref{formula T3(lambda;Nr,Hr)}), and then using the formula (\ref{formula poryadok Hk}), we obtain
\begin{align*}
|T(\alpha;N_\nu,H_\nu)|\ll \frac{H}{N^\frac{n-1}nq^\frac1n}\cdot\frac{qN^\frac{n-1}n}H+q^{\frac12+\varepsilon}\ll q^\frac{n-1}n\ll\\
\ll\left(\frac{H_r}\lnc\right)^\frac{n-1}n\ll\left(\frac{H}{n\mu_k^{1-\frac1n}N^{1-\frac1n}\lnc}\right)^\frac{n-1}n\ll
\frac{H^\frac{n-1}n}{N^\frac{(n-1)^2}{n^2}\lnc^{1-\frac1n}}.
\end{align*}
Substituting the last estimate into (\ref{formula I(M2)-2}), we find
\begin{align*}
&J(\mathfrak{M}_2)\ll\frac{H^{r-1}}{N^{r-\frac rn}}\cdot\frac{N^{n-\frac1n-\frac{n-1}n\varepsilon}}{H^{n-\varepsilon}}\cdot \frac{H^\frac{n-1}n}{N^\frac{(n-1)^2}{n^2}\lnc^{1-\frac1n}}=\frac{H^{r-1}}{N^{r-\frac rn}\lnc^{1-\frac1n}} \left(\frac{N^{\frac{\deg(N)}{\deg(H)}}}H\right)^{\deg(H)},\\
&\deg(N)=n-\frac1n-\frac{n-1}n\varepsilon-\frac{(n-1)^2}{n^2}=n-1+\frac1n-\varepsilon-\frac1{n^2}+\frac1n\varepsilon,\\
&\deg(H)=n-\varepsilon-\frac{n-1}n=n-1+\frac1n-\varepsilon,\\
&\frac{\deg(N)}{\deg(H)}=1-\frac{1-n\varepsilon}{n^3-n^2+n-n^2\varepsilon}=1-\frac1{n^3-n^2+n}+\eta(n),\\
&\eta(n)=\frac{n^2-n}{(n^2-n+1)(n^2-n+1-n\varepsilon)}\varepsilon\le\frac{\varepsilon}{n^2-n+1}.
\end{align*}
From here, taking into account that
$$
1-\frac1{n^3-n^2+n}+\eta(n)\le1-\theta(n)+0,5\varepsilon,
$$
we find
\begin{align}\nonumber
J(\mathfrak{M}_2)&\ll\frac{H^{r-1}}{N^{r-\frac rn}\lnc^{1-\frac1n}}\left(\frac{N^{1-\frac1{n^3-n^2+n}+\eta(n)}}H\right)^{\deg(H)}\ll \\
&\ll\frac{H^{r-1}}{N^{r-\frac rn}\lnc^{1-\frac1n}} \left(\frac{N^{1-\theta(n)+\varepsilon}}H\right)^{\deg(H)}N^{-0,5\varepsilon\deg(H)} \ll\frac{H^{r-1}}{N^{r-\frac rn}\lnc^{r-1}}.
\label{formula I(M2)-3}
\end{align}

{\bf Case 2.} In this case, both conditions of the consequence \ref{Sledst2 Teor ob poved kor trig summi Weyl} of theorem \ref{Teor ob poved kor trig summi Weyl} are satisfied for the sum $T(\alpha;N_\nu,H_\nu)$, that is
$$
\tau=2(n-1)nN_1^{n-2}H_1\ge2(n-1)nN_\nu^{n-2}H_\nu, \qquad\dfrac1{2qN_\nu^{n-1}}<|\lambda |\le\dfrac1{q\tau}.
$$
According to this corollary, conditions $H_r<N_r\lnc^{-1}$ and $q\le H_r\lnc^{-1}$, as well as the relationships (\ref{formula poryadok Nk}) and (\ref{formula poryadok Hk}), we have:
\begin{align*}
|T(\alpha&;N_\nu,H_\nu)|\ll q^\frac{n-1}n\ln q+\min\left(H_\nu q^{-\frac1n},N_\nu^\frac12q^{\frac12-\frac1n}\right)\le\\
&\le q^{1-\frac1n}\ln q+N_\nu^\frac12q^{\frac12-\frac1n}\ll H_r^{1-\frac1n}\lnc^\frac1n+N_r^\frac12H_r^{\frac12-\frac1n}\lnc^{-\frac12+\frac1n}=\\
&=N_r^\frac12H_r^{\frac12-\frac1n}\lnc^{-\frac12+\frac1n}\left(H_r^\frac12N_r^{-\frac12}\lnc^\frac12+1\right)\ll\\
&\ll N_r^\frac12H_r^{\frac12-\frac1n}\lnc^{-\frac12+\frac1n}\ll H^{\frac12-\frac1n}N^{-\frac12+\frac2n-\frac1{n^2}}\lnc^{-\frac12+\frac1n}.
\end{align*}
From here and from (\ref{formula I(M2)-2}), we find
\begin{align*}
J(\mathfrak{M}_2)&\ll\frac{H^{r-1}}{N^{r-\frac rn}}\cdot\frac{N^{n-\frac1n-\frac{n-1}n\varepsilon}}{H^{n-\varepsilon}}\cdot H^{\frac12-\frac1n}N^{-\frac12+\frac2n-\frac1{n^2}}\lnc^{-\frac12+\frac1n}=\\
&=\frac{H^{r-1}}{N^{r-\frac rn}\lnc^{\frac12-\frac1n}} \left(\frac{N^{\frac{\deg(N)}{\deg(H)}}}H\right)^{\deg(H)},\\
\deg(N)&=n-\frac1n-\frac{n-1}n\varepsilon-\frac12+\frac2n-\frac1{n^2}=n-\frac12+\frac1n-\varepsilon-\frac1{n^2}+\frac1n\varepsilon,\\
\deg(H)&=n-\frac12+\frac1n-\varepsilon,\\
\frac{\deg(N)}{\deg(H)}&=1-\frac{1-n\varepsilon}{n^3-0,5n^2+n-n^2\varepsilon}=1-\frac1{n^3-0,5n^2+n}+\eta(n),\\
\eta(n)&=\frac{n^2-0,5n}{(n^2-0,5n+1)(n^2-0,5n+1-n\varepsilon)}\varepsilon%\le\frac{\varepsilon}{n^2-0,5n+1}
 \le\frac{\varepsilon}{n^2-n+1}.
\end{align*}
Hence, keeping in mind that
$$
1-\frac1{n^3-0,5n^2+n}+\eta(n)\le1-\theta(n)+0,5\varepsilon,
$$
we find
\begin{align}\nonumber
J(\mathfrak{M}_2)&\ll\frac{H^{r-1}}{N^{r-\frac rn}\lnc^\frac{n-1}n}\left(\frac{N^{1-\frac1{n^3-0,5n^2+n}+\eta(n)}}H\right)^{\deg(H)}\ll \\
&\ll\frac{H^{r-1}}{N^{r-\frac rn}\lnc^\frac{n-1}n} \left(\frac{N^{1-\theta(n)+\varepsilon}}H\right)^{\deg(H)}N^{-0,5\deg(H)} \ll\frac{H^{r-1}}{N^{r-\frac rn}\lnc^{r-1}}.
\label{formula I(M2)-4}
\end{align}

\subsection{Estimation of the integral $J(\mathfrak{m})$}
Proceeding similarly, as in the case of estimating $J(\mathfrak{M}_2)$, we have
\begin{align}
J(\mathfrak{m})\ll\frac{H^{r-1}}{N^{r-\frac rn}}\cdot\frac{N^{n-\frac1n-\frac{n-1}n\varepsilon}}{H^{n-\varepsilon}}\max_{\alpha\in \mathfrak{m}} |T(\alpha;N_\nu,H_\nu)|,\qquad 1\le\nu\le r. \label{formula I(m)-1}
\end{align}
Let's estimate $T(\alpha;N_\nu,H_\nu)$ for $\alpha$ from the set $\mathfrak{m}$. If $\alpha\in\mathfrak{m}$, then
$$
\alpha=\frac aq+\lambda,\quad (a,q)=1, \quad |\lambda|\le\frac1{q\tau},\quad \frac{H_r}\lnc<q\le\tau,\quad \tau=2n(n-1)N_1^{n-2}H_1.
$$
Using Lemma \ref{Lemma-Weyl-3}, denoting $r=2^n+1$, and then using the relationships
$$
H_\nu\asymp\frac H{N^{1-\frac1n}},\qquad N_\nu\asymp N^\frac1n,\qquad  \frac H{N^{1-\frac1n}\lnc}\ll q\ll\frac H{N^\frac1n},
$$
which are consequences of formulas (\ref{formula poryadok Hk}) and (\ref{formula poryadok Nk}), we have
\begin{align*}
T(\alpha;&N_\nu,H_\nu)\ll H_\nu^{1+\varepsilon}\left(\frac1{H_\nu}+\frac1q+\frac q{H_\nu^n}\right)^\frac2{r-1}\ll\\
&\ll\frac{H^{1+\varepsilon}}{N^{\left(1-\frac1n\right)(1+\varepsilon)}}\left(\frac{N^{1-\frac1n}\lnc}{H}+\frac H{N^\frac1n}\cdot \frac{N^{n-1}}{H^n}\right)^\frac2{r-1}=\\
&=\frac{H^{1-\frac{2n-2}{r-1}+\varepsilon}} {N^{\left(1-\frac1n\right)(1+\varepsilon)-\left(n-1-\frac1n\right)\frac2{r-1}}}\left(\frac{H^{n-2}\lnc}{N^{n-2}}+1\right)^\frac2{r-1}.
\end{align*}
Substituting this estimate into (\ref{formula I(m)-1}), we find
\begin{align}
J(\mathfrak{m})\ll&\frac{H^{r-1}}{N^{r-\frac rn}}\cdot\frac{N^{n-\frac1n-\frac{n-1}n\varepsilon}}{H^{n-\varepsilon}}\cdot  \frac{H^{1-\frac{2n-2}{r-1}+\varepsilon}}{N^{\left(1-\frac1n\right)(1+\varepsilon)-\left(n-1-\frac1n\right)\frac2{r-1}}}=\nonumber\\
 &=\frac{H^{r-1}}{N^{r-\frac rn}}\left(\frac{N^{\frac{\deg(N)}{\deg(H)}}}H\right)^{\deg(H)}, \label{formula I(m)-2}
\end{align}
where
\begin{align*}
 \deg(N) &=\frac1{r-1}\left((r+1)(n-1)-(2r-2)\varepsilon-\frac2n+\frac2n(r-1)\varepsilon\right),\\
\deg(H)&=\frac{rn-r+n-1}{r-1}-2\varepsilon=\frac{(r+1)(n-1)-(2r-2)\varepsilon}{r-1}.
\end{align*}
Using the values of $\deg(N)$ and $\deg(H)$, as well as the notation $r-1=2^n$, we have
\begin{align*}
\frac{\deg(N)}{\deg(H)}&=1-\frac{2-(2r-2)\varepsilon}{(r+1)(n^2-n)-(2r-2)n\varepsilon}=1-\theta(n,r)+\eta(n),\\
\eta(n)&=\frac2{(r+1)(n^2-n)}-\frac{2-(2r-2)\varepsilon}{(r+1)(n^2-n)-(2r-2)n\varepsilon}=\\
&=\frac{\left(1-\frac2{(r+1)(n-1)}\right)(2r-2)\varepsilon}{(r+1)(n^2-n)-(2r-2)n\varepsilon} =\frac{\left(1-\frac2{(2^n+2)(n-1)}\right)2^{n+1}\varepsilon}{(n-1+2^{1-n}(n-1)-2\varepsilon)2^nn}\le\\
&\le\frac{2\varepsilon}{(n-1)n}\le\frac{\varepsilon}2.
\end{align*}
From this estimate and (\ref{formula I(m)-2}), using the relationship $H=N^{1-\theta(n,r)+\varepsilon}$, we find
\begin{align*}
J(\mathfrak{m})&\ll\frac{H^{r-1}}{N^{r-\frac rn}}\left(\frac{N^{1-\theta(n,r)+\frac{\varepsilon}2}}H\right)^{\deg(H)}=\\
&=\frac{H^{r-1}}{N^{r-\frac rn}}\left(\frac{N^{1-\theta(n,r)+\varepsilon}}H\right)^{\deg(H)} N^{-\frac{\deg(H)}2\varepsilon}\ll\frac{H^{r-1}}{N^{r-\frac rn}\lnc^{r-1}}.
\end{align*}
Substituting the obtained estimates for $J(\mathfrak{M}_1)$, $J(\mathfrak{M}2)$, and $J(\mathfrak{m})$ respectively from (\ref{formula I(M1)-3}), (\ref{formula I(M2)-3}), (\ref{formula I(M2)-4}), and (\ref{formula I(m)-2}) into (\ref{formula J_{3,9}(N,H)=I(M1)+I(M2)+I(m)+O(..)}), we obtain the statement of the theorem \ref{TeorAsForWaringPPSl}

\end{document}